%% file: main.tex
\documentclass[12pt]{article}

\include{header}
\newcommand*\mystrut[1]{\vrule width0pt height0pt depth#1\relax}

\newcommand\undermatt[2]{%
  \makebox[0pt][l]{$\smash{\underbrace{\mystrut{3.5ex}\phantom{%
    \begin{matrix}#2\end{matrix}}}_{\text{$#1$}}}$}#2}

\title{A Practical Guide to Optimizing Industrial Thermal Spraying through Comparative Multi-Objective Optimization}

\author{
Wolfgang Rannetbauer\footnote{voestalpine Stahl GmbH, voestalpine-Stra{\ss}e 3, A-4020 Linz, Austria (wolfgang.rannetbauer@voestalpine.at), Corresponding author.} ,
Simon Hubmer\footnote{Johannes Kepler University Linz, Institute of Industrial Mathematics, Altenbergerstra{\ss}e 69, A-4040 Linz, Austria, (simon.hubmer@jku.at)} ,
\\
Carina Hambrock\footnote{voestalpine Stahl GmbH, voestalpine-Stra{\ss}e 3, A-4020 Linz, Austria (carina.hambrock@voestalpine.com)} ,
Ronny Ramlau\footnote{Johannes Kepler University Linz, Institute of Industrial Mathematics, Altenbergerstra{\ss}e 69, A-4040 Linz, Austria, (ronny.ramlau@jku.at)} \footnote{Johann Radon Institute for Computational and Applied Mathematics, Altenbergerstra{\ss}e 69, A-4040 Linz, Austria, (ronny.ramlau@ricam.oeaw.ac.at)}
}

\begin{document}

\maketitle
\begin{abstract}

Achieving both high quality and cost-efficiency are two critical yet often conflicting objectives in manufacturing and maintenance processes. Quality standards vary depending on the specific application, while cost-effectiveness remains a constant priority. These competing objectives lead to multi-objective optimization problems, where algorithms are employed to identify Pareto-optimal solutions—compromise points which provide decision-makers with feasible parameter settings. The successful application of such optimization algorithms relies on the ability to model the underlying physical system, which is typically complex, through either physical or data-driven approaches, and to represent it mathematically. This paper applies three multi-objective optimization algorithms to determine optimal process parameters for high-velocity oxygen fuel (HVOF) thermal spraying. Their ability to enhance coating performance while maintaining process efficiency is systematically evaluated, considering practical constraints and industrial feasibility. Practical validation trials are conducted to verify the approximate theoretical solutions generated by the algorithms, ensuring their applicability and reliability in real-world scenarios. By exploring the performance of these diverse algorithms in an industrial setting, this study offers insights into their practical applicability, guiding both researchers and practitioners in enhancing process efficiency and product quality in the coating industry.

\end{abstract}

\section{Introduction}
In the modern era of advanced manufacturing and industrial processes, the search for optimal performance is central to technological innovation and operational excellence \cite{hughes2022perspectives}. Achieving this level of optimization often requires more than just refining individual components or processes; it demands a deep understanding of the complex systems at hand, typically achieved by developing models which capture the sophisticated dependencies within these systems \cite{jabbari2018multiphysics}. This comprehensive understanding not only facilitates the realization of economic advantages, such as enhanced cost efficiency and waste reduction, but also contributes to the improvement of product quality, including attributes such as durability and resistance. 

A particularly challenging and relevant area where optimization plays a critical role is thermal spraying, a flexible and widely used technique for applying coatings to surfaces \cite{pierlot2008design, fauchais2014industrial}. Traditionally, manufacturers of coating materials have developed sets of standard parameters, which serve as preliminary guidelines or approximations for the application of various materials \cite{venkatachalapathy2023guiding}. However, these sets of parameters are frequently based on generalized assumptions and may not exist for all materials. Moreover, the effectiveness of these parameters is influenced by the specific hardware configuration employed—factors such as the type of spraying system, the capacity of available compressed air, gas pressures, and even regional variations in gas composition can lead to significant discrepancies in outcomes \cite{fauchais2014thermal}. Consequently, these sets of standard parameters often prove inadequate, particularly when the coatings are intended for use under atypical or extreme conditions, as is frequently encountered in industries such as steel manufacturing. This deficiency highlights a significant optimization challenge: the need to precisely tailor process parameters to achieve the desired coating properties—such as hardness, adhesion, and surface finish—while also considering economic factors like material usage and process efficiency \cite{tillmann2010desirability}.

In thermal spraying, the optimization problem becomes even more complex due to the multitude of interdependent variables—such as spray distance, particle velocity, and substrate temperature—that can influence the final coating quality. The development of accurate mathematical models which represent these complex interactions is essential for the successful application of optimization algorithms \cite{rannetbauer2024enhancing,rannetbauer2024predictive,rannetbauer2024virtual,rannetbauer2024leveraging}. These models, whether derived from physical principles or data-driven approaches, provide a simplified yet powerful representation of the thermal spraying process, enabling the identification of optimal parameter settings \cite{calafiore2014optimization}. 

Building upon this foundation, quality characteristics and economic aspects are treated as objectives which the model seeks to optimize. In optimization, these objectives can be singular, focusing on a single goal such as minimizing cost, or multiple, where several goals, such as balancing cost with quality, must be optimized simultaneously. The latter, known as multi-objective optimization, is particularly challenging and requires algorithms capable of navigating trade-offs between competing objectives \cite{deb2016multi}. In the context of thermal spraying, this involve optimizing parameters to achieve a coating that meets stringent performance criteria while also minimizing material waste and energy consumption.

This work investigates three established multi-objective optimization algorithms to address the specific challenges encountered in thermal spray coating technology. While these algorithms have been widely employed across various industrial optimization problems, their adaptation to the thermal spraying domain offers unique opportunities for generating practical insights. By tailoring these methods to the specific requirements of thermal spraying, this study focuses on comparing their performance and suitability for optimizing critical coating properties under realistic operational constraints. The findings aim to provide a practical framework for industry practitioners and researchers engaged in advanced manufacturing processes.

This paper first provides a brief overview of the physical principles of thermal spraying, followed by a precise definition of its primary optimization challenges in Section~\ref{sect_background}. Section~\ref{sect_mathMod} introduces the mathematical modelling framework, employing generalized linear models (GLMs) to represent the relationships between input parameters and coating quality attributes. Section \ref{sect_optimTheory} outlines optimization theory, focusing on multi-objective approaches and key concepts such as dominance and the Pareto front. Section~\ref{sect_optimMethods} examines the theoretical basis and methodology of three selected multi-objective optimization algorithms, which are subsequently applied to practical thermal spraying problems in Section \ref{sect_application}. Validation results assessing the effectiveness of these approaches are presented in Section \ref{sect_validation}, followed by a summary of practical implications and recommendations in Section \ref{sect_conclusion}.

\section{Physical Background and Problem Definition}\label{sect_background}
Thermal spray coating technology is an indispensable process in a range of industrial applications, providing protective and functional coatings which improve the performance and durability of components subjected to demanding environments \cite{davis2004handbook, fauchais2014thermal}. In steel manufacturing, thermal spray coatings play a critical role in extending the service life of equipment by providing resistance to wear, corrosion, and high-temperature degradation. As a result, thermal spray technology is employed across various stages of steel production, including steel making, continuous casting, rolling, and other high-wear processes which operate under extreme conditions \cite{singh2023applications}. For example, in 2008, 20\% of thermal spray applications in China were used in the steel industry \cite{fukumoto2008current}. While more recent data for other regions are unavailable, this illustrates the significance of thermal spray technology in steel manufacturing.

The flexibility of thermal spray technology lies in its ability to apply a wide variety of materials—including metals, ceramics, and composites—onto substrates without significantly altering the underlying base material properties \cite{amin2016review}. Several distinct thermal spray processes have been developed over the years, each with its specific advantages depending on the application requirements. These processes include plasma spraying, arc spraying, flame spraying, and cold spraying, among others \cite{mathesius2009praxis}. Each of these methods relies on the use of an energy source to melt or semi-melt coating materials, which are then propelled onto a substrate at high speeds to form a compact, adherent coating. The suitability of a particular thermal spray process is largely determined by the characteristics of the coating material—such as its type, particle size, and chemical composition—as well as the intended application of the coated component \cite{ang2014review}. The processes differ primarily in the levels of thermal and kinetic energy they can impart, making certain methods more appropriate for specific materials and performance requirements.

\subsection{High-Velocity Oxy-Fuel Spraying}
Among these methods, HVOF (High-Velocity Oxy-Fuel) spraying, an advancement of flame spraying, has emerged as one of the most effective techniques for producing high-quality coatings with superior bond strength and low porosity \cite{davis2004handbook}. Due to its capability to produce coatings which resist extreme environments, HVOF is particularly appropriate for diverse applications in the steel industry. 

\begin{figure}[htb!]
    \centering
    \includegraphics[width=0.90\textwidth]{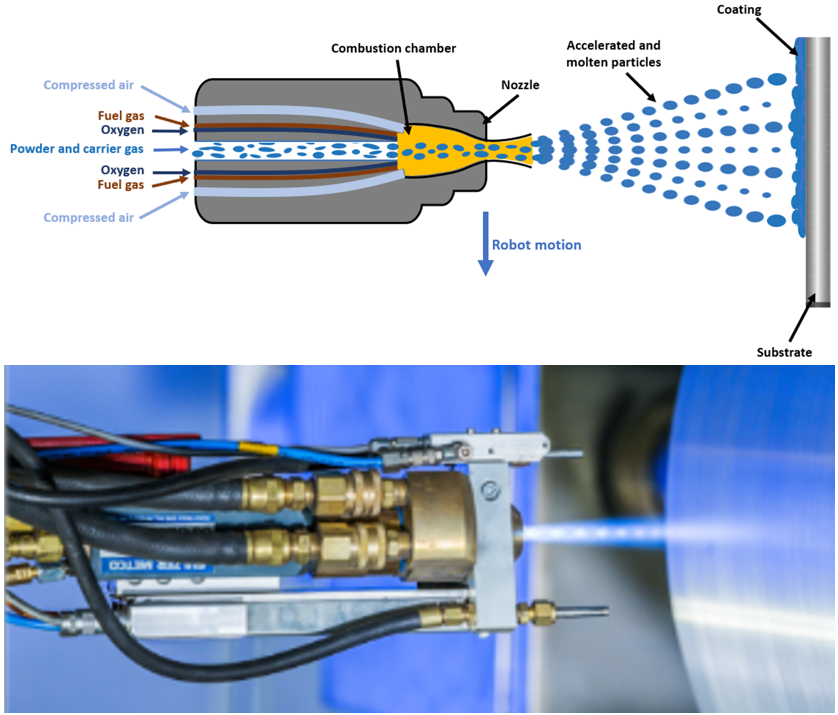}
    \caption{Schematic and real-world representation of a gas-fueled HVOF spraying process. Top: Schematic of the HVOF spray gun, showing the axial injection of powder into the combustion chamber, where oxygen and fuel gases react. The heated particles are accelerated through an expansion nozzle and deposited onto the substrate.
    Bottom: Photograph of the HVOF process, taken from \cite{TSM2}, displaying the characteristic shock diamonds, indicating supersonic velocities and high kinetic energy.}
    \label{fig:hvof_general}
\end{figure}

The HVOF process, illustrated in Figure \ref{fig:hvof_general}, is characterized by the combustion of a fuel gas—typically hydrogen, propane, natural gas, or liquid fuels like kerosene—with oxygen in a high-pressure combustion chamber. This exothermic reaction generates extremely high temperatures and pressures, forcing the combustion gases through an expansion nozzle to accelerate the gas flow to supersonic velocities. Compressed air is additionally injected to further increase the speed of the gas stream, while the coating material, as powder, is introduced axially (gas-fuel) into the system via an inert carrier gas. As the particles pass through the combustion chamber and the expansion nozzle, they are heated and strongly accelerated before exiting the spray gun. Upon impact with the substrate, these molten particles flatten, adhere to the pre-treated surface, and solidify, forming a dense, cohesive coating. The mechanical bond formed by the high-velocity impact is responsible for the superior adhesion and durability of HVOF coatings \cite{fauchais2014thermal}.

In the bottom panel of Figure \ref{fig:hvof_general}, the characteristic shock diamonds, visible within the spray plume, are evidence of supersonic flow, highlighting the high kinetic energy transmitted to the particles during this process. The HVOF gun is typically mounted on the end of a robotic arm, allowing precise and controlled movement, which is important for achieving uniform coatings across various geometries. Additionally, it is cooled using water and/or air to maintain operational stability and prevent overheating.

Further essential information, particularly on the technical key aspects relevant to mathematical modelling and optimization in thermal spraying, is summarized in \cite{rannetbauer2024predictive}. For a more comprehensive understanding of the fundamentals of thermal spraying and detailed insights into the HVOF process, readers are encouraged to consult \cite{davis2004handbook, fauchais2014thermal}, which offer in-depth discussions on coating preparation, process principles, parameters, and technological advancements in coating techniques.
 
While the HVOF process involves a multitude of input parameters, research has shown that only a subset of these parameters has a significant influence on both coating quality and process efficiency~\cite{becker2021artificial, liu2019prediction, rannetbauer2024enhancing, ribu2022experimental, tillmann2022statistical}. These input parameters, including the powder feed rate, the stand-off distance, the fuel-to-oxygen ratio, the coating velocity, and the total gas flow rate (which in some studies is further broken down into the individual flow rates of fuel gas (or kerosene), oxygen, and compressed air), directly impact key coating properties such as porosity, hardness, and wear resistance, as well as process performance metrics like deposition efficiency. Consequently, adjusting these parameters determines not only the quality, durability, and functional performance of the coating in industrial applications but also the cost-effectiveness of the coating process itself.

However, each input parameter is subject to physical boundaries dictated by equipment limitations, material properties, and safety considerations. For instance, exceeding optimal powder feed rates can lead to incomplete particle melting, resulting in higher coating roughness and poorer adhesion \cite{rannetbauer2024predictive}, while excessively high stand-off distances may cause a significant drop in particle velocity, negatively impacting coating density and hardness \cite{tillmann2022statistical}. Similarly, deviations in the fuel-to-oxygen ratio can alter flame temperature and particle acceleration, further influencing the microstructure and mechanical properties of the coating.

Thus, it is essential to operate within the defined physical boundaries of the input parameters to not risk poorer coating performance, process inefficiency, or equipment wear. Depending on the coating material used, specific parameter limits may vary, reflecting material-dependent factors like particle size, melting point, and thermal conductivity. Table \ref{tab:limits} outlines general and material-specific boundaries for the five key input parameters, mentioned above, in HVOF spraying. The material-specific boundaries are subject to a tungsten carbide powder (WC-10Co-4Cr) with a grain size in the range of -45 + 15 µm, supplied by Oerlikon Metco \cite{metco2023material}. This type of coating material is frequently applied in the steel industry due to its excellent wear resistant properties.

\begin{table}[h] 
    \resizebox{\textwidth}{!}{
    \begin{tabular}{llcccc}  
    \toprule
    \multicolumn{2}{c}{Input parameter} & \multicolumn{2}{c}{General boundaries} & \multicolumn{2}{c}{WC-10Co-4Cr boundaries}\\
    \cmidrule(r){1-2} \cmidrule(r){3-4} \cmidrule(r){5-6}
    \multicolumn{1}{l}{Name}    & \multicolumn{1}{l}{Unit} & \multicolumn{1}{c}{Lower} & \multicolumn{1}{c}{Upper} & \multicolumn{1}{c}{Lower} & \multicolumn{1}{c}{Upper}\\
    \midrule
    Powder feed rate (PFR)      & g/min    & 20     & 150 & 45     & 75\\
    Stand-off distance (SOD)    & mm    & 170     & 320 & 200     & 260\\
    Fuel-to-oxygen ratio ($\lambda$)  & -    & 0.50     & 1.20 & 0.84     & 1.04\\
    Coating velocity (CV)       &  m/min    & \undermatt{\text{highly equipment-specific}}{ \ \ \ \ \ 50 \ \ \ \ \ &  \ \ 150  \ \  &} 75 & 125 \\
    Total gas flow (TGF)      &  nl/m    & 530 & 830 & 615 & 751 \\
    \\
    \\
    \bottomrule
    \end{tabular}
    }
    \caption{Physical boundaries of input parameters in HVOF spraying, with general limits shown in the center and material-specific boundaries for WC-10Co-4Cr powder presented on the right.}
    \label{tab:limits}
\end{table}

While the general boundaries presented in Table \ref{tab:limits} provide a useful reference for HVOF process parameters, it is important to note that these limits are highly equipment-specific and not universally applicable. The general values are based on the operational capabilities of an Oerlikon Metco thermal spraying system, specifically the DJ 2700 gas-fuel HVOF system with a water-cooled gun assembly, a Metco 9MPE-DJ powder feeder, and a standard gas supply management system using propane as fuel gas. Moreover, the fuel-to-oxygen ratio is particularly sensitive to the coating material, since different materials require specific flame temperatures and combustion environments to achieve optimal particle melting and deposition. Variations in melting point, thermal conductivity, and particle behavior under high heat influence the choice of the fuel-to-oxygen ratio for each material system. These limits, therefore, serve as a baseline for processes conducted on similar setups, while deviations may be required for other equipment configurations or materials.

Optimizing HVOF process parameters is essential for achieving application-specific coating properties, particularly in industrial applications such as steel manufacturing, where coatings must resist harsh operating conditions. The challenge is to determine the optimal combination of input parameters—such as powder feed rate, stand-off distance, fuel-to-oxygen ratio, coating velocity, and total gas flow—that yields coatings with the (multiple) desired properties while remaining within the physical boundaries dictated by equipment and material constraints. Specifically, the coating must be optimized for its mechanical and chemical properties \cite{fauchais2014thermal, rannetbauer2024predictive}, while also maximizing process efficiency to ensure durability, operational reliability, and cost-effectiveness.

\section{Mathematical Modelling with Generalized Linear Models}\label{sect_mathMod}
In prior research, predictive models for critical HVOF coating properties were developed using generalized linear models (GLMs) based on a systematic design of experiments, specifically a central composite design (CCD)  \cite{rannetbauer2024predictive}, involving the five HVOF input parameters outlined in Table \ref{tab:limits}. This paper builds upon those insights to determine application-specific optimal HVOF input parameters. 

GLMs represent an extension of traditional linear regression techniques, allowing for the modelling of response variables which are assumed to follow a distribution from the exponential family, such as Gaussian, binomial, or gamma distributions \cite{nelder1972generalized}. Unlike ordinary linear regression, GLMs model the relationship between a dependent variable and one or more independent variables through a link function which connects the linear predictor to the mean of the response variable. The general form of a GLM is given by:
    \begin{align}
        g(\mu_i) = h_{\beta}(\boldsymbol{x}_i),
    \end{align}
where $\mu_i = \mathbb{E}(y_i|\boldsymbol{x}_i)$ represents the conditional mean of the response variable $y_i$ for observation $i$. Collectively, the vector of conditional means across all observations is $\boldsymbol{\mu} = (\mu_i)_{i=1}^n$, with the response variable $\boldsymbol{y} = (y_i)_{i=1}^n$ representing the coating properties of interest. Here, $\mathbf{x}_i = (1, x_{1}^{i},x_{2}^{i},\dots,x_{k}^{i})_{i=1}^n$ is the vector of predictor variables. The function $h_{\beta}(\boldsymbol{x}_i)$ includes linear, quadratic, and interaction terms of the predictor variables:
\begin{align}
h_\beta(\boldsymbol{x}_i) = \beta_0 + \sum_{j=1}^{k} \beta_j {x}^{i}_j + \sum_{j=1}^{k} \sum_{m=j}^{k} \beta_{jm} {x}^{i}_j {x}^{i}_m,
\end{align}
where $\beta_0$ is the intercept, $\beta_j$ are the coefficients for the linear terms, $\beta_{jj}$ are the coefficients for the quadratic terms, and $\beta_{jm}$ are the coefficients for the interaction terms. The link function $g(\cdot)$ serves to relate the mean of the response variable $\mu_i$ to the predictor function $h_\beta(\boldsymbol{x}_i)$.

In the case of thermal spray coating properties, the gamma distribution was adopted due to its suitability for handling variables that are strictly positive and exhibit a non-symmetrical (conditional) distribution, as observed in all evaluated coating characteristics \cite{rannetbauer2024predictive}. For gamma regression, the natural choice for the link function is the logarithmic function, ensuring that the predicted values remain positive \cite{fahrmeir2013generalized}. The gamma regression model can be expressed as: 
    \begin{align} \label{eqGamma}
        \log(\mu_i) = h_{\beta}(\boldsymbol{x}_i).
    \end{align}
The explicit formulas for the final models, which accurately describe eight distinct characteristics of HVOF coatings, are provided in Appendix \ref{appendixA}. These characteristics include particle in-flight properties (velocity, temperature), process performance properties (deposition rate, deposition efficiency), and coating quality properties (thickness, roughness, hardness, porosity). Particle in-flight properties are critical, since they influence the resultant microstructure, affecting whether certain phases are developed and whether the coating adopts an amorphous or crystalline structure. Process performance metrics quantify the efficiency of the deposition process. Finally, coating quality properties, such as thickness, roughness, hardness, and porosity, are direct indicators of the coating's suitability for its intended application, impacting mechanical properties like wear resistance, adhesion, and durability.

To ensure comparability and avoid scale effects, all input parameters were normalized according to the CCD using factors, which accounts for the smaller coefficient magnitudes in Appendix \ref{appendixA}. 
For the HVOF coating properties, the likelihood function was derived from the gamma distribution, and the estimation was performed using maximum likelihood estimation \cite{rannetbauer2024predictive}. This estimation procedure quantifies the regression coefficients that are most probable given the observed data $(y_i,\boldsymbol{x}_i)_{i=1}^n$, considering the assumed conditional distribution of the response variable $\boldsymbol{y}$ and the selected link function $g(\cdot)$.

To refine the models and ensure their suitability for predicting coating properties, an iterative model selection process was employed. Model selection was guided by the Akaike Information Criterion (AIC), which balances the model fit and complexity by penalizing excessive model parameters \cite{akaike1974new}. This approach led to reduced models which provided a concise yet accurate representation of the relationship between HVOF input parameters and coating properties, as evidenced by high adjusted $R^2$ values and minimal cross-validation errors \cite{rannetbauer2024enhancing, rannetbauer2024predictive}. The resulting final models serve as the basis for the optimization of application-specific HVOF input parameters. 

In the context of industrial-scale operations in the steel industry, optimizing these properties ensures that both the functional requirements of the coating and the efficiency of the process are met. By focusing on these parameters, the optimization framework captures the essential characteristics required to achieve application-specific performance while maintaining process efficiency.

\section{Multi-objective Optimization Theory}\label{sect_optimTheory}

\begin{figure}[t!]
    \centering
    \includegraphics[width=1\textwidth]{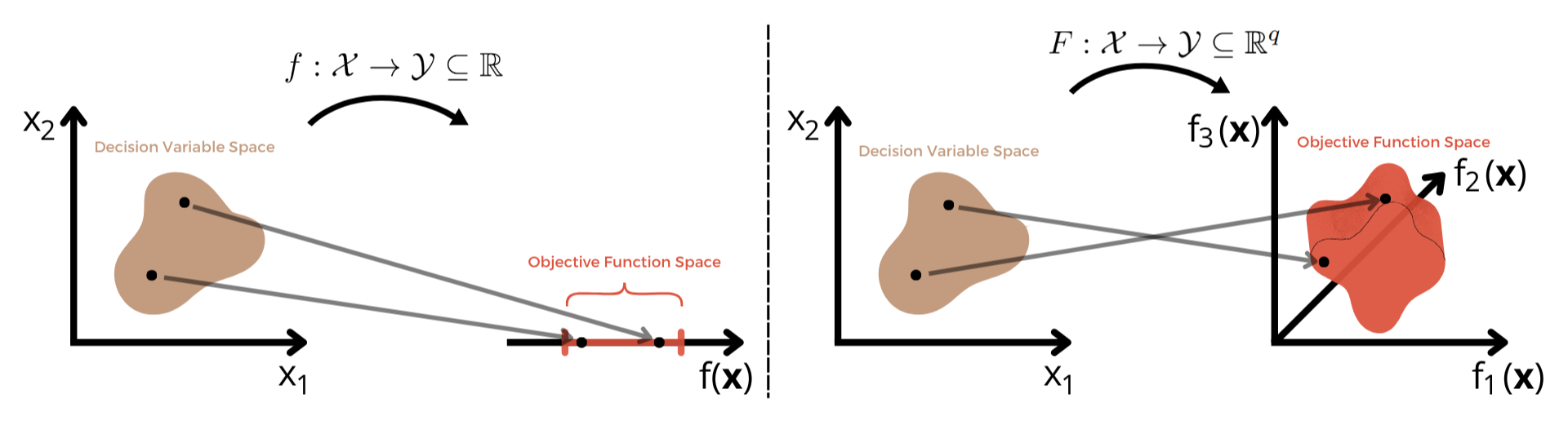}
    \caption{Comparison of single-objective and multi-objective optimization. The left panel shows the mapping from the decision variable space \(\mathcal{X}\) to a one-dimensional objective function space \(\mathcal{Y} \subseteq \mathbb{R}\). The right panel illustrates the mapping to a higher-dimensional objective space \(\mathcal{Y} \subseteq \mathbb{R}^{q}\). Here $\mathcal{Y}$ is three-dimensional.}
    \label{fig:multi_overrview}
\end{figure}

Multi-objective optimization involves the simultaneous optimization of multiple conflicting objectives, often requiring trade-offs between them \cite{gunantara2018review}. These problems typically include constraints that any feasible solution, including the optimal ones, \linebreak must satisfy. Formally, the objective functions are collected into a vector \linebreak $F(\boldsymbol{x}) = \left[f_1(\boldsymbol{x}), f_2(\boldsymbol{x}), \dots, f_{q}(\boldsymbol{x})\right]^\top$. The resulting solution landscape is inherently broader and more complex, often yielding a set of trade-off solutions rather than a single optimal point \cite{gunantara2018review}. These problems operate within two essential Euclidean spaces \cite{chiandussi2012comparison}:
\begin{itemize}
    \item The $n$-dimensional decision variable space $\mathcal{X} \subseteq \mathbb{R}^n$, where each axis corresponds to a component of the decision vector $\boldsymbol{x}$.
    \item The $q$-dimensional objective space $\mathcal{Y} \subseteq \mathbb{R}^{q}$, where each axis represents one of the objective functions $f_l(\boldsymbol{x})$.
\end{itemize}
This distinction is visually represented in Figure \ref{fig:multi_overrview}, where the mapping from the decision space to the objective space is shown for both single-objective and multi-objective cases. The formal definition of a multi-objective optimization problem is as follows \cite{deb2016multi}:
\begin{definition}[Multi-Objective Optimization Problem]\label{multi-optim}
    A multi-objective optimization problem is defined as the task of minimizing a vector of objective functions $F(\boldsymbol{x}) = \left[f_1(\boldsymbol{x}), f_2(\boldsymbol{x}), \dots, f_{q}(\boldsymbol{x})\right]^\top$ subject to a set of inequality constraints $g_i(\boldsymbol{x})$ and equality constraints $h_j(\boldsymbol{x})$, formulated as:
    \begin{align} \label{multi-prob}
    \begin{split}
    \text{minimize} \quad & F(\boldsymbol{x}) = \left[f_1(\boldsymbol{x}), f_2(\boldsymbol{x}), \dots, f_{q}(\boldsymbol{x})\right]^\top \\
    \text{subject to} \quad & g_i(\boldsymbol{x}) \leq 0, \quad \forall \ i \in { 1,\dots, m} \\
    & h_j(\boldsymbol{x}) = 0, \quad \forall \ j \in { 1,\dots, p}
    \end{split}
    \end{align}
    where $\boldsymbol{x} \in \mathcal{X} \subseteq \mathbb{R}^n$ is an $n$-dimensional vector of decision variables defined over the input space $\mathcal{X}$. The objective functions $f_l(\boldsymbol{x})$, for $l=1,\dots, q$, each map the decision variable vector $\boldsymbol{x}$ to a scalar value, with the vector $F(\boldsymbol{x})$ representing the collective output in the $q$-dimensional objective function space $\mathcal{Y} \subseteq \mathbb{R}^{q}$.
\end{definition}

\begin{remark}
    As in single-objective optimization, maximization in multi-objective optimization can be reformulated as a minimization problem by negating the objective function, i.e., multiplying $f_l(\boldsymbol{x})$ by -1, in accordance with the duality principle \cite{rao1984optimization}. The same principle can be applied to inequality constraints. If a constraint is expressed as a ``less-than-or-equal-to" condition, it can be transformed into an equivalent ``greater-than-or-equal-to" form by multiplying the constraint function by -1 \cite{deb2016multi}. This allows both types of constraints to be handled consistently within the optimization framework.

      The input space $\mathcal{X}$, sometimes referred to as the domain of the optimization problem \eqref{multi-prob}, is defined as the intersection of the domains of all objective and constraint functions, specifically \cite{boyd2004convex}:
    \begin{align} \label{feasible}
        \mathcal{X} = \bigcap_{i=1}^m \text{dom}(g_i) \cap \bigcap_{j=1}^p \text{dom}(h_j).
    \end{align}
    A solution $\boldsymbol{x} \in \mathcal{X}$ is considered feasible if it satisfies all constraints, meaning $g_i(\boldsymbol{x}) \leq 0$ for $i = 1,\dots,m$, and $h_j(\boldsymbol{x}) = 0$ for $j = 1,\dots,p$. Feasibility ensures that only those solutions meeting the predefined constraints are considered, which is fundamental to the optimization process \cite{deb2016multi}. Multi-objective optimization is often referred to as vector optimization, reflecting the simultaneous optimization of a vector of objectives rather than a single scalar objective \cite{deb2016multi}. Thus, trade-offs between conflicting objectives are typically required.
\end{remark}

\begin{definition}[Ideal Objective Vector]
The ideal objective vector, denoted by $\boldsymbol{p}^\star$, represents the theoretical optimum in a multi-objective optimization problem, where each component of $\boldsymbol{p}^\star$ corresponds to the best possible value for an individual objective function, considered in isolation. Formally, the $l$-th component of the ideal objective vector, $p^\star_l$, is defined as:
\begin{align}
    p_l^\star = \inf \big \{ f_l(\boldsymbol{x}) \ \mid \ \boldsymbol{x} \in \mathcal{X}, \ g_i(\boldsymbol{x}) \leq 0, \ i = 1, \dots, m, \ h_j(\boldsymbol{x}) = 0, \ j = 1, \dots, p \big \} 
\end{align}
for each objective $l = 1, \dots, q$. The ideal objective vector is then given by:
\begin{align}
    \boldsymbol{p}^\star = \left( p_1^\star, p_2^\star, \dots, p_{q}^\star \right)^\top.
\end{align}

In general, $\boldsymbol{p}^\star$ does not correspond to any feasible solution in the decision space $\mathcal{X}$, as the optimal values $p_l^\star$ for each objective function are typically achieved at different points $\boldsymbol{x}$. The ideal objective vector represents an utopian point where all objectives would achieve their optimal values simultaneously, which is generally unattainable unless the objectives are non-conflicting. Nevertheless, it serves as a fundamental benchmark in multi-objective optimization, guiding the search for optimal trade-offs and enabling solution comparison. Consequently, solutions closer to $\boldsymbol{p}^\star$ are generally preferred \cite{deb2016multi}.
\end{definition}

\subsection{Concept of domination}
In multi-objective optimization, the concept of domination is essential for comparing solutions across multiple, often conflicting objectives. It provides a systematic method for determining whether one solution is superior to another, serving as the basis for many optimization algorithms.

\begin{definition}[Weak Dominance] \label{pareto_dom}
Consider two candidates, $\boldsymbol{x}^{(1)}$ and $\boldsymbol{x}^{(2)}$ within the decision space $\mathcal{X}$. The candidate $\boldsymbol{x}^{(1)}$ is considered to weakly dominate $\boldsymbol{x}^{(2)}$, denoted by $\boldsymbol{x}^{(1)} \preceq \boldsymbol{x}^{(2)}$, if and only if the following two conditions are satisfied \cite{nagy2020multi}:
\begin{enumerate}
\item For all objectives $l \in \{1, 2, \dots, q\}$, the objective function value corresponding to $\boldsymbol{x}^{(1)}$ is no worse than that corresponding to $\boldsymbol{x}^{(2)}$, i.e.,
\begin{align}\label{dom1}
    f_l(\boldsymbol{x}^{(1)}) \leq f_l(\boldsymbol{x}^{(2)}), \quad \forall \ l \in \{1, 2, \dots, q\}.
\end{align}
\item There exists at least one objective $l \in \{1, 2, \dots, q\}$ for which the objective function value corresponding to $\boldsymbol{x}^{(1)}$ is strictly better than that corresponding to $\boldsymbol{x}^{(2)}$, i.e.,
\begin{align}\label{dom2}
    \exists \ l \in \{1, 2, \dots, q\}: \ f_l(\boldsymbol{x}^{(1)}) < f_l(\boldsymbol{x}^{(2)}).
\end{align}
\end{enumerate}
If both conditions are met, $\boldsymbol{x}^{(1)}$ is said to weakly dominate $\boldsymbol{x}^{(2)}$, indicating that $\boldsymbol{x}^{(2)}$ is suboptimal when compared to $\boldsymbol{x}^{(1)}$. This definition assumes a minimization context, such as \eqref{multi-prob}, where lower objective function values are preferred. In cases involving maximization or a mix of objectives, the criteria for domination would need to be appropriately adjusted.
\end{definition}

\begin{definition}[Strong Dominance] \label{strong_dom}
Consider two candidates, $\boldsymbol{x}^{(1)}$ and $\boldsymbol{x}^{(2)}$, within the decision space $\mathcal{X}$. The candidate $\boldsymbol{x}^{(1)}$ is said to strongly dominate $\boldsymbol{x}^{(2)}$, denoted by $\boldsymbol{x}^{(1)} \prec \boldsymbol{x}^{(2)}$, if and only if:
\begin{align}\label{strong_dom_cond}
f_l(\boldsymbol{x}^{(1)}) < f_l(\boldsymbol{x}^{(2)}), \quad \forall \ l \in {1, 2, \dots, q}.
\end{align}
\end{definition}
This means that $\boldsymbol{x}^{(1)}$ performs strictly better than $\boldsymbol{x}^{(2)}$ across all objectives, without exception. In other words, $\boldsymbol{x}^{(1)}$ is preferable to $\boldsymbol{x}^{(2)}$ in every respect.

\subsection{Concept of Pareto-Optimality} \label{pareto-optimality}
In multi-objective optimization, the concept of Pareto-optimality provides a crucial framework for identifying solutions which achieve the best possible trade-offs among conflicting objectives. 

\begin{definition}[Pareto-Optimality] \label{pareto_opt}
A candidate $\boldsymbol{x}^\star \in \mathcal{X}$ is considered Pareto-optimal if there is no other solution $\boldsymbol{x} \in \mathcal{X}$ which weakly dominates $\boldsymbol{x}^\star$. \end{definition}

This definition implies that no other feasible solution can improve an objective without causing at least one other objective to become worse.

\begin{definition}[Pareto Optimal Set] \label{pareto_set}
The Pareto Optimal Set, denoted by $\mathcal{P}^\star \subseteq \mathcal{X}$, is the set of all Pareto-optimal solutions within the decision space $\mathcal{X}$. Formally, the Pareto Optimal Set is defined as:
\begin{align}
\mathcal{P}^\star = \left\{ \boldsymbol{x}^\star \in \mathcal{X} \mid \nexists \ \boldsymbol{x} \in \mathcal{X} \text{ such that } \boldsymbol{x} \preceq \boldsymbol{x}^\star \right\}.
\end{align}
\end{definition}

\begin{figure}[htb!]
    \centering
    \includegraphics[width=0.85\textwidth]{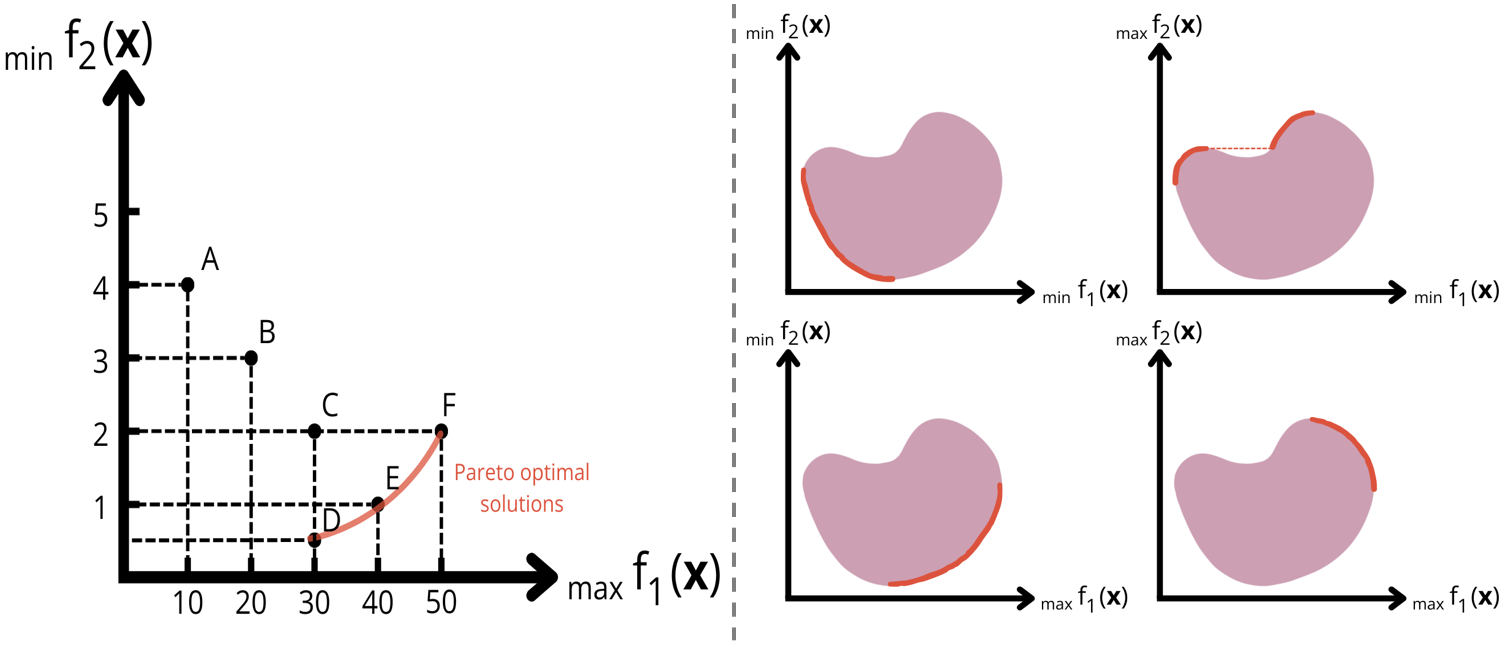}
    \caption{Left: Illustration of Pareto domination in a two-objective optimization problem. The figure plots six candidates, labeled A through F, in a coordinate system where $f_1$ is maximized and $f_2$ is minimized. Candidates D, E, and F are identified as non-dominated, representing Pareto optimal solutions. The figure is adapted from \cite{deb2016multi}. Right: Visualization of Pareto-Optimal Sets in a Two-Objective Optimization Problem for Different Objective Orientations  (Adapted from \cite{deb2016multi}). Each subfigure illustrates the Pareto-optimal set for distinct scenarios where the objectives $f_1$ and $f_2$ are either minimized or maximized.}
    \label{fig:domination_concept}
\end{figure}

\begin{example}[Dominance and Pareto Optimality in a Two-Objective Optimization Problem] \label{ex_dominance_pareto}\ In a two-objective optimization problem with $f_1(\boldsymbol{x})$ to be maximized and $f_2(\boldsymbol{x})$ to be minimized, six candidates (A through F) are depicted in Figure \ref{fig:domination_concept} (left). These candidates illustrate the trade-offs between conflicting objectives and the concepts of dominance and Pareto optimality.
Candidate B dominates A, since $f_1({B}) > f_1({A})$ and $f_2({B}) < f_2({A})$, fulfilling the conditions for both weak and strong dominance. Similarly, D weakly dominates C, since $f_2({D}) < f_2({C})$ while $f_1({D}) = f_1({C})$, but D does not strongly dominate C due to equality in $f_1$. This distinction highlights that strong dominance implies weak dominance, but the reverse is not true.
Candidates D, E, and F are non-dominated, meaning no improvements in one objective are possible without compromising the other. These candidates form the Pareto Optimal Set, $\mathcal{P}^\star = \{\boldsymbol{x}_D, \boldsymbol{x}_E, \boldsymbol{x}_F\}$, representing the best trade-offs achievable for $f_1$ and $f_2$. Here, $\boldsymbol{x}_D$, $\boldsymbol{x}_E$, and $\boldsymbol{x}_F$ are vectors in the $n$-dimensional decision space, each corresponding to a specific configuration of decision variables.
\end{example}

\begin{remark}
The set of objective vectors corresponding to all candidates in the Pareto Optimal Set $\mathcal{P}^\star$ forms what is known as the Pareto Front, denoted by $\mathcal{F}^\star$. Mathematically, the Pareto Front is defined as:
\begin{align}
\mathcal{F}^\star = \left\{ F(\boldsymbol{x}) \mid \boldsymbol{x} \in \mathcal{P}^\star \right\}.
\end{align}

\end{remark}

In the two-objective optimization problem discussed in Example~\ref{ex_dominance_pareto}, the Pareto Front $\mathcal{F}^\star$ is visualized in the objective space, where each point represents the objective function values of a Pareto-optimal candidate from $\mathcal{P}^\star$. This is also illustrated in the left plot of Figure~\ref{fig:domination_concept}, where the projections of candidates D, E, and F form the Pareto Front. The right plot further depicts Pareto-optimality under different optimization directions for $f_1$ and $f_2$: minimizing both objectives (top-left), minimizing $f_1$ and maximizing $f_2$ (top-right), maximizing $f_1$ and minimizing $f_2$ (bottom-left), and maximizing both (bottom-right). Depending on the optimization direction, the Pareto-optimal set lies along a distinct boundary segment of the objective space.

Building on the concepts of Pareto-optimality and the Pareto Front, the objectives in multi-objective optimization are twofold:

\begin{enumerate}
    \item Proximity to the Pareto Optimal Set: It is essential to find solutions which closely approximate the Pareto Optimal Set, $\mathcal{P}^\star$. This ensures solutions near the optimal trade-offs, similar to the goal in single-objective optimization.
    \item Diversity of Solutions: The second goal, unique to multi-objective optimization, is to ensure a diverse representation of solutions within the Pareto-optimal region. Diversity in this context means that the solutions should be well-spread across the Pareto Front $\mathcal{F}^\star$ in the objective space $\mathcal{Y}$. A diverse set of solutions provides a comprehensive view of the trade-offs between objectives, enabling better decision-making \cite{deb1999multi}.
\end{enumerate}
In summary, while convergence to a solution set near the Pareto-optimal front is critical for ensuring near-optimality, achieving a diverse set of solutions is crucial for exploring the full spectrum of trade-offs between objectives.

\section{Optimization Methods}\label{sect_optimMethods}
In industrial applications, such as the optimization of complex processes like thermal spraying, determining appropriate process parameters is critical to achieving desired quality characteristics. Multi-objective optimization methods provide a structured framework for addressing such problems by systematically balancing competing objectives. A widely accepted classification of these methods is based on the role of the decision-maker and the availability of preference information, distinguishing four principal categories \cite{miettinen1999nonlinear, deb2016multi}:

\begin{enumerate}
    \item \textbf{\textit{No-preference methods:}} These methods operate independently of the decision-maker, identifying a compromise solution from the Pareto Optimal Set $\mathcal{P}^\star$ without requiring preference information. Such methods are suitable when no specific preferences are available or when the involvement of the decision-maker is limited.
    \item \textbf{\textit{A priori methods:}} Here, the decision-maker specifies preferences before optimization begins, guiding the formulation of the objective function. The quality of the solution depends on the clarity and accuracy of these initial preferences \cite{miettinen1999nonlinear}.
    \item \textbf{\textit{A posteriori methods:}} These approaches generate a representative set of Pareto-optimal solutions $\mathcal{P}^\star$, allowing the decision-maker to evaluate and select a solution post-optimization. While flexible, these methods can become computationally demanding, especially for problems with many objectives or complex solution spaces \cite{deb2016multi}.
    \item \textbf{\textit{Interactive methods:}} Also known as \textit{Progressive Preference Articulation} \cite{adra2007comparative, chiandussi2012comparison, rostami2017progressive}, these methods involve iterative feedback from the decision-maker, refining the optimization process to align solutions with evolving preferences.
\end{enumerate}

The methods explored in this work are selected to represent different categories of multi-objective optimization techniques, ensuring a comprehensive approach to the problem. These methods, along with their respective classifications, are summarized in Table~\ref{tab:method_classification}.

\begin{table}[h!] 
\centering 
\resizebox{\textwidth}{!}{
\begin{tabular}{lcccc} \hline \textbf{Optimization Method} & \textbf{No-preference} & \textbf{A priori} & \textbf{A posteriori}& \textbf{Interactive} \\ 
\hline 
Weighted Sum Method with SQP & & \checkmark & \checkmark & \\
Desirability Approach & & \checkmark & & \\ 
Non-dominated Sorting Genetic Algorithm II & & & \checkmark & \\ 
\hline 
\end{tabular}
}
\caption{Classification of optimization methods based on decision-maker involvement and preference information.} 
\label{tab:method_classification} 
\end{table}

These methods are chosen for their practical relevance and suitability for solving multi-objective optimization problems in industrial applications. No-preference methods were not considered, as selecting appropriate HVOF parameters requires explicitly incorporating performance criteria relevant to coating quality and process efficiency. Moreover, interactive methods were omitted due to the impracticality of iterative decision-maker input in an industrial environment, where optimization must be performed efficiently without continuous manual adjustments. Detailed descriptions of the selected methods, including their theoretical foundations and implementation, are provided in the following sections.

\subsection{Weighted Sum Method with SQP Optimization}
The Weighted Sum Method, introduced in \cite{gass1955computational}, is among the most widely employed techniques in multi-objective optimization due to its simplicity and compatibility with established single-objective optimization algorithms \cite{marler2010weighted}. This method reformulates a multi-objective problem by aggregating the objectives into a single scalar function: 
\begin{align} \label{weighted-sum}
    \begin{split}
    \underset{\boldsymbol{x} \in \mathcal{X}}{\text{minimize}} & \quad \sum_{l=1}^{q} w_l f_l(\boldsymbol{x}),\\
    \text{subject to} & \quad g_i(\boldsymbol{x}) \leq 0, \quad \forall \ i \in \{ 1,\dots, m\} \\
    & \quad h_j(\boldsymbol{x}) = 0, \quad \forall \ j \in \{ 1,\dots, p\}
    \end{split}
\end{align}
where $\boldsymbol{w} = [w_1, \dots, w_{q}]^\top$ is a weight vector satisfying $w_l \geq 0$ for all $l \in \{1, \dots, q\}$ and $\sum_{l=1}^{q} w_l = 1$. The optimization is performed over the input space $\mathcal{X}$, incorporating both inequality constraints $g_i(\boldsymbol{x}) \leq 0$ and equality constraints $h_j(\boldsymbol{x}) = 0$ as defined in~\eqref{feasible}.

This method can be employed in two primary ways, distinguished by the timing of preference articulation by the decision-maker. In an a priori approach, the decision-maker specifies the weight vector $\boldsymbol{w}$ prior to optimization. This ensures that the computed solution reflects the predefined priorities \cite{miettinen1999nonlinear}. Alternatively, in an a posteriori framework, a series of weight vectors is systematically applied to generate multiple solutions along the Pareto front, from which the decision-maker can later select \cite{deb2016multi}.

Despite its widespread adoption, the Weighted Sum Method is not without limitations. It is known to be ineffective in capturing solutions in non-convex regions of the Pareto front, since the scalarized objective inherently assumes linear combinations of objectives, which cannot represent such regions \cite{das1997closer,nakayama1995aspiration}. Additionally, even if weights are uniformly distributed, the corresponding Pareto-optimal solutions may not be evenly spaced along the front, especially in problems with complex geometries or highly nonlinear objectives \cite{das1997closer,deb2016multi}. 

Given the nature of the scalarized objective in the Weighted Sum Method, various single-objective optimization techniques can be applied to solve problem~\eqref{weighted-sum}. Among these, Sequential Quadratic Programming (SQP) is particularly well-suited due to its numerical stability and well-established convergence properties in constrained nonlinear optimization \cite{nocedal1999numerical}. SQP iteratively solves a sequence of quadratic programming (QP) subproblems which approximate the original problem by linearizing constraints and using a second-order approximation of the objective function. The method ensures convergence under suitable regularity conditions and is particularly effective for problems with smooth nonlinear objectives and constraints \cite{bonnans2006numerical}. At each iteration, SQP minimizes a quadratic approximation of the Lagrangian, defined as:

\begin{align} \label{lagrangian}
    \begin{split}
    \mathcal{L}(\boldsymbol{x},\boldsymbol{\lambda},\boldsymbol{\nu}) = f(\boldsymbol{x}) + \sum_{i=1}^{m}\lambda_ig_i(\boldsymbol{x}) + \sum_{j=1}^{p}\nu_jh_j(\boldsymbol{x}),
    \end{split}
\end{align}
where $f(\boldsymbol{x})$ denotes the objective function, and $\boldsymbol{\lambda} = (\lambda_i)_{i=1}^m$ and $\boldsymbol{\nu} = (\nu_j)_{j=1}^p$ are the Lagrange multipliers associated with the inequality and equality constraints, respectively. To solve the optimization problem, SQP approximates the Lagrangian $ \mathcal{L}(\boldsymbol{x},\boldsymbol{\lambda},\boldsymbol{\nu})$ using a second-order Taylor expansion around the current iterate $\boldsymbol{x}_{\tau}$, leading to the following QP subproblem:

\begin{align} \label{taylor}
    \begin{split}
        \underset{\boldsymbol{d}}{\text{minimize}} \quad & f(\boldsymbol{x}_{\tau}) + \nabla f(\boldsymbol{x}_{\tau})^\top \boldsymbol{d} + \frac{1}{2} \boldsymbol{d}^\top \nabla^2 \mathcal{L}(\boldsymbol{x}_{\tau}, \boldsymbol{\lambda}_{\tau}, \boldsymbol{\nu}_{\tau})\boldsymbol{d}, \\
        \text{subject to} \quad & g_i(\boldsymbol{x}_{\tau}) + \nabla g_i(\boldsymbol{x}_{\tau})^\top \boldsymbol{d} \leq 0, \quad  \forall \ i \in \{ 1,\dots, m\}, \\
        & h_j(\boldsymbol{x}_{\tau}) + \nabla h_j(\boldsymbol{x}_{\tau})^\top \boldsymbol{d} = 0, \quad  \forall \ j \in \{ 1,\dots, p\}.
    \end{split}
\end{align}
Here, $\nabla^2 \mathcal{L}(\boldsymbol{x}_{\tau}, \boldsymbol{\lambda}_{\tau}, \boldsymbol{\nu}_{\tau})$ is the Hessian matrix of the Lagrangian, and $\nabla f(\boldsymbol{x}_{\tau})$, $\nabla g_i(\boldsymbol{x}_{\tau})$, and $\nabla h_j(\boldsymbol{x}_{\tau})$ are the gradients of the objective and constraints at $\boldsymbol{x}_{\tau}$. The Hessian of the Lagrangian, $\nabla^2 \mathcal{L}(\boldsymbol{x}_{\tau}, \boldsymbol{\lambda}_{\tau}, \boldsymbol{\nu}_{\tau})$, is composed of second derivatives of the objective function and constraints. To reduce computational complexity, quasi-Newton methods such as BFGS are commonly used to approximate the Hessian matrix $\nabla^2 \mathcal{L}(\boldsymbol{x}_{\tau}, \boldsymbol{\lambda}_{\tau}, \boldsymbol{\nu}_{\tau})$ by a positive definite symmetric matrix $\boldsymbol{H}_{\tau}$, ensuring positive definiteness while leveraging curvature information from gradients \cite{pan2016numerical, broyden1970convergence, fletcher1970new, goldfarb1970family, shanno1970conditioning}. The iterative update is expressed as \cite{nocedal1999numerical}:
\begin{sloppypar}
\begin{align} \label{BFGS}
    \boldsymbol{H}_{{\tau}+1} = \boldsymbol{H}_{\tau} -\frac{\boldsymbol{H}_{\tau} \boldsymbol{s}_{\tau} \boldsymbol{s}_{\tau} ^\top \boldsymbol{H}_{\tau}}{\boldsymbol{s}_{\tau} ^\top \boldsymbol{H}_{\tau} \boldsymbol{s}_{\tau}} + \frac{\boldsymbol{y}_{\tau} \boldsymbol{y}_{\tau}^\top}{\boldsymbol{y}_{\tau}^\top \boldsymbol{s}_{\tau}},
\end{align}
where $\boldsymbol{s}_{\tau} = \boldsymbol{x}_{{\tau}+1} - \boldsymbol{x}_{\tau}$ is the difference between consecutive iterates (= step direction), and ${\boldsymbol{y}_{\tau} = \nabla \mathcal{L}(\boldsymbol{x}_{{\tau}+1}) - \nabla \mathcal{L}(\boldsymbol{x}_{\tau})}$ is the difference in gradients of the Lagrangian. 
\end{sloppypar}

The resulting iterative procedure involves solving the QP subproblem, updating the decision variables $\boldsymbol{x}_{\tau}$, and refining the Hessian approximation $\boldsymbol{H}_{\tau}$. The algorithm continues until the convergence criteria are met, such as when changes in the decision variables fall below a specified tolerance.

To ensure the applicability of SQP in the context of the weighted sum approach, it is necessary to verify that the objective function in \eqref{weighted-sum} is differentiable and, ideally, convex. Differentiability follows directly from the properties of the individual objective functions $f_l(\boldsymbol{x})\overset{\eqref{eqGamma}}{=}\exp{(h_\beta(\boldsymbol{x}))}$, which, as established in Section~\ref{sect_mathMod}, correspond to gamma regression models with a logarithmic link function. The gradient of $f_l(\boldsymbol{x})$ is given by:

\begin{align*}
    \nabla f_l(\boldsymbol{x}) = \exp{(h_\beta(\boldsymbol{x})) \nabla h_\beta(\boldsymbol{x})},
\end{align*}
where
\begin{align*}
    \nabla h_\beta(\boldsymbol{x}) = \Big( \frac{\partial h_\beta(\boldsymbol{x})}{\partial x_1}, \dots, \frac{\partial h_\beta(\boldsymbol{x})}{\partial x_k} \Big)
\end{align*}
Since these terms are well-defined and continuous, differentiability of the scalarized function is ensured. 

Regarding convexity, the Hessian of $f_l(\boldsymbol{x})$ takes the form: 
\begin{align*}
    H = \nabla^2 f_l(\boldsymbol{x}) = \exp(h_\beta(\boldsymbol{x})) (\boldsymbol{\beta} \boldsymbol{\beta}^T + B(\boldsymbol{x})),     
\end{align*}
where \( B(\boldsymbol{x}) \) is the Hessian of \( h_\beta(\boldsymbol{x}) \), which depends on the quadratic and interaction terms. The entry $B_{jm}(\boldsymbol{x})$ of the matrix $B(\boldsymbol{x})$ is given by:
\begin{align*}
    B_{jm}(\boldsymbol{x}) = \frac{\partial^2}{\partial x_j \partial x_m} (\beta_0 + \sum_{j=1}^{k} \beta_j x_j \sum_{j=1}^{k} \sum_{m=j}^{k} \beta_{jm} x_j x_m).
\end{align*}
Evaluating this gives:
\begin{align*}
    B_{jm}(\boldsymbol{x}) =
    \begin{cases}
        \beta_{jm}, & \text{if } j \neq m, \\
        2\beta_{jj}, & \text{if } j = m.
    \end{cases}
\end{align*}

Since $\exp(h_\beta(\boldsymbol{x})) > 0$, the convexity of $ f_l(\boldsymbol{x})$ is dictated by the definiteness of $ \boldsymbol{\beta} \boldsymbol{\beta}^T + B(\boldsymbol{x}) $. If interaction effects are absent, $h_\beta(\boldsymbol{x})$ reduces to a purely quadratic function in the input parameters, where each term contributes a positive definite component to \( B(\boldsymbol{x}) \), ensuring positive definiteness of $ \boldsymbol{\beta} \boldsymbol{\beta}^T + B(\boldsymbol{x}) $. However, when interaction terms are present, \( B(\boldsymbol{x}) \) may introduce indefinite components, potentially violating convexity. 

Since all model functions are explicitly known, it is possible to compute $ \boldsymbol{\beta} \boldsymbol{\beta}^T + B(\boldsymbol{x}) $ for each combination of gamma regression models and determine positive definiteness analytically. In cases where convexity is not guaranteed, quasi-Newton methods such as BFGS within SQP remain applicable, since they approximate the Hessian while ensuring numerical stability even in non-convex settings.

As an alternative to SQP for solving the scalarized, box-constrained optimization problems, the Spectral Projected Gradient (SPG) method presents an efficient option~\cite{birgin2000nonmonotone, birgin2014spectral}. The method combines the classical projected gradient approach with a spectral step-length strategy that incorporates second-order information without computing the Hessian matrix~\cite{barzilai1988two}. A key feature of this approach is a non-monotone line search, which permits temporary increases in the objective function. This behavior can be advantageous for navigating complex optimization landscapes more effectively.

A comparative analysis conducted for this work revealed that SPG leads to solutions that are, for all practical purposes, identical to those obtained by SQP. While SPG was found to be computationally somewhat faster, this advantage is of limited practical relevance in the context of this study, given that the absolute optimization time is negligible compared to the overall industrial process time. Nevertheless, the method remains a valuable alternative, especially for potential future work on problems with a higher number of parameters where its scalability would be beneficial.

\subsection{Desirability Approach}
The desirability function approach, while less sophisticated from a mathematical perspective, remains one of the most extensively utilized methods in industry for optimizing multi-response processes \cite{karande2013applications}. Its popularity is largely attributed to its simplicity and intuitive framework, making it a practical choice for industrial applications, particularly in materials science. As an a priori method, it requires the specification of weights, bounds, and other parameters before optimization, thereby encoding the decision-maker's preferences directly into the model. This study compares the desirability approach with mathematical optimization techniques and evolutionary algorithms to assess their respective effectiveness in multi-objective optimization.

Originally introduced by Harrington \cite{harrington1965desirability} and later refined by Derringer and Suich~\cite{derringer1980simultaneous}, the desirability approach transforms multiple conflicting objectives into a unified optimization framework. For each objective \(f_j(\boldsymbol{x})\) in the multi-objective optimization problem, the desirability function assigns a value \(d_j \in [0, 1]\), where \(d_j = 0\) represents an entirely unacceptable outcome, and \(d_j = 1\) corresponds to the ideal or desired performance \cite{derringer1980simultaneous, he2020robust}. The transformation depends on whether the goal is to minimize or maximize a specific value of the objective. The formulations are as follows:

\subsubsection{One-Sided Transformation}
For objectives requiring either minimization or maximization, the desirability functions are defined as:
\paragraph{Minimization Objective}
\begin{align}
    d_j(f_j(\boldsymbol{x})) = 
    \begin{cases} 
    0, & f_j(\boldsymbol{x}) \geq U_j, \\
    \left(\frac{U_j - f_j(\boldsymbol{x})}{U_j - L_j}\right)^r, & L_j \leq f_j(\boldsymbol{x}) < U_j, \\
    1, & f_j(\boldsymbol{x}) \leq L_j,
    \end{cases}
\end{align}
where \(L_j\) and \(U_j\) are the lower and upper bounds for the objective function \(f_j(\boldsymbol{x})\), respectively, and the exponent \(r > 0\) is a shape parameter controlling the curvature of the desirability function.

\paragraph{Maximization Objective}
\begin{align}
    d_j(f_j(\boldsymbol{x})) = 
    \begin{cases} 
    0, & f_j(\boldsymbol{x}) \leq L_j, \\
    \left(\frac{f_j(\boldsymbol{x}) - L_j}{U_j - L_j}\right)^r, & L_j < f_j(\boldsymbol{x}) \leq U_j, \\
    1, & f_j(\boldsymbol{x}) \geq U_j,
    \end{cases}
\end{align}
where \(L_j\), \(U_j\), and \(r\) are as defined above. 

\begin{remark}
    The choice of the shape parameter \(r\) depends on technical, economic, or other practical considerations \cite{costa2011desirability}. For \(r = 1\), the desirability function exhibits a linear increase or decrease towards the upper or lower limit. Considering the minimization of an objective, \(r < 1\) results in a convex function $d_j(f_j(\boldsymbol{x}))$, favoring values closer to \(L_j\). Conversely, \(r > 1\) results in a concave function. This parameter provides flexibility to adapt the desirability function to the specific needs of the optimization problem at hand.
\end{remark}

\subsubsection{Overall Desirability Function}

The overall desirability function \(D\) provides a scalar measure to assess the optimization solution across multiple objectives. It is defined as the geometric mean of the individual desirability values \(d_j\), ensuring that \(D\) is sensitive to the performance of all objectives. The mathematical formulation is given by:
\begin{align}
    D(\boldsymbol{x}) = \left(\prod_{j=1}^{q} d_j(f_j(\boldsymbol{x}))^{\rho}_j\right)^{\frac{1}{\sum_{j=1}^{q} \rho}_j},
\end{align}
where \(d_j \in [0, 1]\) represents the individual desirability for the \(j\)-th objective, \(\rho_j \geq 0\) denotes the weight assigned to the \(j\)-th objective, and \(\sum_{j=1}^{q} \rho_j > 0\). The weights \(\rho_j\) allow for flexibility in assigning importance to different objectives, reflecting specific technical, economic, or practical priorities.

This formulation ensures that the overall desirability \(D\) lies within the interval \([0, 1]\). A value of \(D = 0\) occurs when at least one objective is entirely unacceptable (\(d_j = 0\)), whereas \(D = 1\) corresponds to the ideal case where all objectives achieve maximum desirability. The use of the geometric mean, rather than an arithmetic mean, ensures that the measure is sensitive to poor performance in any single objective, preventing the overall assessment from being heavily influenced by favorable outcomes in others. 

The geometric mean formulation also exhibits the property of scale invariance, which is critical for aggregating dimensionless measures of desirability. Additionally, the inclusion of weights allows the decision-maker to emphasize specific objectives, providing a mechanism to adapt the optimization process to varying priorities or constraints. For equal weights (\(\rho_j = \frac{1}{j}\) for all \(j\)), the formulation reduces to a simple geometric mean, reflecting an unbiased aggregation of individual desirability values.

Despite the lack of differentiability in the individual desirability functions, the method achieves robust optimization through direct search techniques, which are well-suited for non-differentiable optimization landscapes \cite{kolda2003optimization}. The desirability approach has demonstrated its efficiency in optimizing multiple responses in thermal spraying processes, including kerosene-fueled HVOF systems \cite{tillmann2010desirability}, and is expected to exhibit similar performance in gas-fueled systems.

Its widespread adoption in industrial and materials science applications stems from its flexibility in accommodating diverse quality characteristics and its ability to prioritize specific performance criteria \cite{costa2011desirability}. Unlike traditional multi-objective optimization methods, which typically minimize or maximize aggregated objective functions as expressed in \eqref{multi-prob}, the desirability approach offers the additional capability of targeting specific response values through two-sided transformations \cite{derringer1980simultaneous}. In this study, however, only one-sided transformations are employed to ensure a consistent basis for comparison with other optimization methods.

\subsection{Non-dominated Sorting Genetic Algorithm II}
The Non-dominated Sorting Genetic Algorithm II (NSGA-II), introduced in \cite{deb2002fast}, represents an advanced evolutionary optimization algorithm designed to efficiently address multi-objective problems. As a refinement of the original NSGA framework \cite{srinivas1994muiltiobjective}, NSGA-II integrates the principles of genetic algorithms (GA) with a robust non-dominated sorting approach \cite{simon2013evolutionary}. A key improvement over its predecessor is the incorporation of the crowding distance operator, which effectively preserves solution diversity—a critical objective in multi-objective optimization, as outlined in Section~\ref{pareto-optimality}—and addresses a known limitation of the original NSGA \cite{deb2002fast}. Although various modified versions of NSGA-II have been developed \cite{deb2013evolutionary, asefi2014hybrid}, this study employs the original algorithm due to its demonstrated robustness \cite{dang2023analysing}, conceptual simplicity, and proven applicability across a wide range of optimization scenarios \cite{verma2021comprehensive}.

The NSGA-II algorithm operates by iteratively evolving a population of candidate solutions through genetic operators such as selection, crossover, and mutation. Its core mechanism is structured around two primary objectives: the identification of non-dominated solutions and the preservation of diversity within the Pareto front. 

\subsubsection{Core Components of NSGA-II}
The NSGA-II framework is built upon four fundamental principles: non-dominated sorting, the elite-preserving operator, the crowding distance metric, and the selection mechanism \cite{deb2002fast}.

\begin{figure}[htb!]
    \centering
    \includegraphics[width=0.85\textwidth]{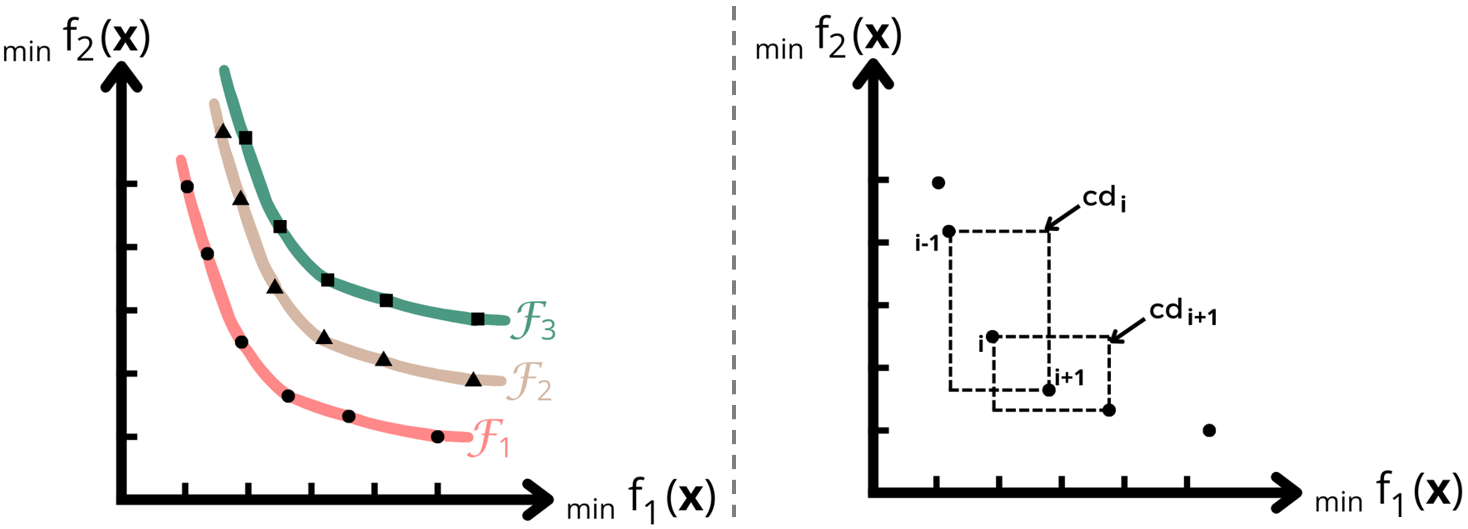}
    \caption{Left: Illustration of the non-dominated sorting procedure for a bi-objective minimization problem, displaying three Pareto fronts: $\mathcal{F}_1$ (rank 1, $\bullet$), $\mathcal{F}_2$ (rank 2, $\blacktriangle$), and $\mathcal{F}_3$ (rank 3, $\blacksquare$). The classification of individuals into hierarchical fronts demonstrates the ranking process based on Pareto dominance. Right: Visualization of crowding distance computation in a bi-objective minimization problem. The crowding distances $\operatorname{cd}_i$ and $\operatorname{cd}_{i+1}$ are depicted using dotted cuboids, representing the differences in objective values between adjacent solutions in the objective space.}
    \label{fig:NSGA_theory}
\end{figure}

\begin{enumerate}
    \item Non-dominated sorting: The initial population is randomly generated within predefined box constraints, ensuring feasibility with respect to variable bounds. Non-dominated sorting classifies the population into hierarchical fronts based on Pareto dominance (Definition~\ref{pareto_dom}). First, all non-dominated candidates are identified and assigned rank 1, forming the first Pareto front, $\mathcal{F}_1$. These candidates are then removed, and the process iterates for the remaining individuals. The second front, $\mathcal{F}_2$, consists of candidates that are non-dominated after excluding $\mathcal{F}_1$, and so forth. This classification continues until all individuals are assigned to a front, as illustrated in the left panel of Figure~\ref{fig:NSGA_theory} for a bi-objective minimization problem.

    \item Elite-preserving operator: NSGA-II integrates an elitism mechanism to ensure that the highest-quality candidates are retained across generations. This strategy saves non-dominated candidates by directly transferring them to the subsequent population, thereby preventing their premature loss during the evolutionary process. Specifically, non-dominated individuals from the current generation are preserved and carried forward until they are dominated by superior candidates in future iterations. This approach ensures steady progression toward an improved Pareto front.
    
    \item Crowding distance: To maintain solution diversity, NSGA-II employs the crowding distance metric, which quantifies the relative spacing of candidates in the objective space. Candidates in less crowded regions are favored during selection, promoting a well-distributed Pareto front.

    The crowding distance \(\operatorname{cd}_i\) of an individual \(i\) is computed as the average normalized difference in objective function values between its nearest neighbors:  
    \begin{align}
        \operatorname{cd}_i = \sum_{l=1}^q \frac{f_l^{i+1} - f_l^{i-1}}{f_l^{\text{max}} - f_l^{\text{min}}},
    \end{align}
    where \(q\) is the number of objective functions, \(f_l^{i}\) denotes the \(l^{\text{th}}\) objective value of candidate \(i\), and \(f_l^{\text{max}}\) and \(f_l^{\text{min}}\) are the maximum and minimum values of the \(l^{\text{th}}\) objective, respectively. 

    The population is sorted for each objective in ascending order before computing crowding distances \cite{deb2002fast}. Boundary solutions, which lack neighbors on one side, are assigned an infinite crowding distance (\(\operatorname{cd}_1 = \operatorname{cd}_{|\mathcal{I}|} = \infty\)) to ensure their retention. Here, \(|\mathcal{I}|\) represents the total number of solutions within the non-dominated set \(\mathcal{I}\). This approach preserves both boundary solutions and those in sparsely populated regions. The right panel of Figure~\ref{fig:NSGA_theory} illustrates the crowding-distance computation for a bi-objective minimization problem.  
    
    \item Selection operator: The selection of individuals for the next generation in NSGA-II is performed using a crowded binary tournament selection operator, which evaluates individuals based on rank and crowding distance. The decision rule follows:  
    \begin{itemize}
    \item The solution with the superior (lower) rank is selected if the competing individuals have different ranks.  
    \item If the individuals share the same rank, the one with the greater crowding distance is chosen to enhance diversity. 
    \end{itemize}
    This selection mechanism ensures a balance between convergence toward the Pareto-optimal front and the preservation of diversity within the solution set.
\end{enumerate}

\subsubsection{Algorithmic Framework of NSGA-II}
The NSGA-II algorithm begins by initializing a random population $P_0$ of size $N$. A secondary population $Q_0$ is then generated by applying binary tournament selection, followed by crossover and mutation operations on $P_0$ \cite{blickle2000tournament}. These genetic operators introduce variation, enabling exploration of the solution space and preventing premature convergence.

\begin{figure}[htb!]
    \centering
    \includegraphics[width=0.85\textwidth]{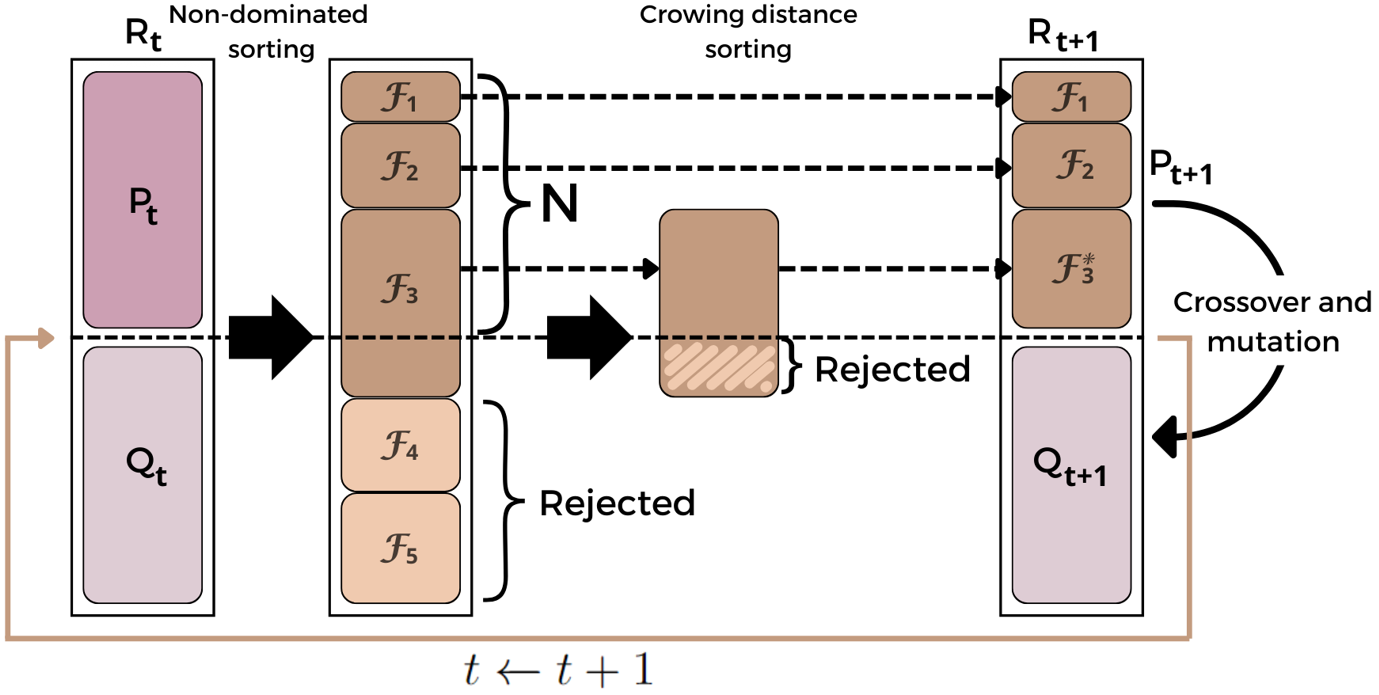}
    \caption{Schematic representation of the NSGA-II algorithmic framework (Adapted from \cite{deb2002fast}).  At each generation $t$, the parent and offspring populations are merged to form a combined population $R_t$, which undergoes non-dominated sorting to classify solutions into Pareto fronts $\mathcal{F}_1, \mathcal{F}_2, \dots$. Solutions are then selected iteratively, prioritizing lower-ranked fronts, until the population size $N$ is reached. To maintain diversity, the final front $\mathcal{F}_j$, here $\mathcal{F}_3$, is sorted by crowding distance, ensuring the selection of well-distributed solutions. The resulting population $P_{t+1}$ is subjected to crossover and mutation operations to produce $Q_{t+1}$, iterating until convergence or a predefined stopping criterion is satisfied.}
    \label{fig:NSGA_procedure}
\end{figure}

The iterative procedure of NSGA-II, outlined in Figure~\ref{fig:NSGA_procedure}, proceeds as follows:
At each generation $t$, the current population $P_t$ is combined with the offspring population $Q_t$ to form a combined population $R_t = P_t \ \cup \ Q_t$. Non-dominated sorting is applied to $R_t$, ranking the solutions into hierarchical fronts $\mathcal{F}_1, \mathcal{F}_2, \dots$, based on their dominance relationships. The first front $\mathcal{F}_1$, containing the highest-ranked non-dominated solutions, is prioritized for the next generation's population $P_{t+1}$. If $\mathcal{F}_1$ contains fewer than $N$ solutions, all its members are included in $P_{t+1}$, and additional members are selected sequentially from subsequent fronts ($\mathcal{F}_2, \mathcal{F}_3, \dots$) until $P_{t+1}$ reaches size $N$.

If including an entire front $\mathcal{F}_j$ would exceed the population size $N$, selection within $\mathcal{F}_j$ is performed based on crowding distance, favoring solutions in less crowded regions of the objective space (i.e., those with higher crowding distance values).

This ensures a balance between convergence and diversity. Once $P_{t+1}$ is finalized, it is subject to selection, crossover, and mutation to produce the next offspring population $Q_{t+1}$. This process iterates until the termination criterion—typically a predefined number of generations—is met.

\section{Application of Optimization Methods to Multi-Objective Optimization Problems in Thermal Spray Coating Technology}\label{sect_application}

Building on the theoretical foundations and optimization methodologies outlined in the preceding sections, this section presents a case study on optimizing HVOF process parameters for industrial coating applications.

Based on the technical setup detailed in Section~\ref{sect_background}, three optimization problems of particular relevance to industrial HVOF thermal spray coating are identified. These problems pertain to the application of tungsten carbide cobalt chrome (WC-10Co-4Cr) powder, a material known for its exceptional abrasion resistance and wear protection. Since coating hardness is the primary determinant of these performance characteristics, it constitutes the principal optimization objective. However, in certain industrial applications, particularly within the steel industry, additional coating properties are critical, requiring a broader optimization approach beyond hardness alone.

The optimization problems are formulated using gamma regression models, as described in Section~\ref{sect_mathMod} and detailed in Appendix~\ref{appendixA}. These models provide accurate predictions of coating properties based on HVOF process parameters, as demonstrated in \cite{rannetbauer2024predictive}. The optimization constraints are defined by material-specific parameter boundaries, summarized in Table~\ref{tab:limits}. These boundaries establish feasible ranges for key process variables, including powder feed rate, stand-off distance, fuel-to-oxygen ratio, coating velocity, and total gas flow, thereby ensuring practical and physically realizable solutions.

The optimization methods are implemented within the \texttt{R} statistical computing environment \cite{R}, leveraging specialized packages tailored to the methodologies under consideration. Specifically, the \texttt{desirability2} \cite{desirability2} package is employed for desirability-based optimization, while \texttt{NLopt} \cite{NLopt} and \texttt{nsga2R} \cite{nsga2R} are utilized for gradient-based and evolutionary multi-objective optimization, respectively.

Beyond box constraints, no additional equality or inequality constraints are explicitly introduced, since the optimization aims to identify feasible solutions within the predefined parameter space. Consequently, the problem formulation reduces to a multi-objective optimization task constrained by physical boundaries, facilitating a systematic exploration of trade-offs among competing objectives such as surface hardness, coating roughness, and deposition efficiency. The following subsections provide a detailed formulation of the identified optimization problems and present the results obtained using the proposed methods.

\subsection{Problem I: Optimization of Surface Hardness and Deposition Efficiency}

The first optimization problem considers the dual objectives of maximizing surface hardness and deposition efficiency, reflecting the general purpose of using tungsten carbide cobalt chrome (WC-10Co-4Cr) coatings. High hardness is critical for ensuring abrasion resistance and wear protection while optimizing deposition efficiency minimizes material waste and enhances process cost-effectiveness. Let \( f_1(\mathbf{x}) \) denote the surface hardness and \( f_2(\mathbf{x}) \) the deposition efficiency, where \( \mathbf{x} = (x_1, x_2, x_3, x_4, x_5) \) represents the process parameters. The optimization problem is formulated as:

\[
\text{Minimize } F(\boldsymbol{x}) = \left[-f_1(\mathbf{x}), -f_2(\mathbf{x})\right]^\top ,
\]
subject to the physical constraints on the process parameters:

\[
x_{\text{Lower}} \leq x \leq x_{\text{Upper}},
\]
where \( x_{\text{Lower}} \) and \( x_{\text{Upper}} \) denote the material-specific parameter boundaries, as detailed in Table~\ref{tab:limits}.

To solve this bi-objective optimization problem, the three optimization methods, described in Section~\ref{sect_optimMethods} are applied: the Weighted Sum Method with SQP optimization, the desirability approach, and the NSGA-II algorithm. Each method requires specific hyperparameter settings to ensure proper implementation. For the Weighted Sum Method, weight vectors are systematically applied in an a posteriori manner to generate a Pareto front approximation. To improve numerical stability and avoid local optima, the optimization is initialized using multiple starting values $\boldsymbol{x}_0$ for the decision variables, sampled from within the feasible parameter space. Convergence is determined based on the relative tolerance for changes in decision variables:

\[
\frac{|\Delta x_i|}{|x_i|} < \texttt{xtol}, \quad \text{where } \texttt{xtol} = 10^{-8}.
\]
Here, \(|\Delta x_i|\) denotes the absolute change in the value of the \(i\)-th decision variable \(x_i\) between successive iterations, and \(|x_i|\) represents the absolute value of the current value of \(x_i\). 

In the desirability approach, desirability functions are defined for both objectives, characterized by lower and upper bounds ($L_j$ and $U_j$, respectively) and a shape parameter ($r_j$) controlling the desirability curve. Objective weights ($\rho_j$) are adjusted to reflect their relative importance. The specific parameter values are detailed in Table~\ref{tab:hyperparameters}.

For the NSGA-II algorithm, a population size of (\(N=300\)) and a maximum of 1000 generations are chosen to balance exploration of the solution space and computational efficiency. 

\begin{table}[ht]
\centering
\resizebox{\textwidth}{!}{
\begin{tabular}{llll}
\hline
\textbf{Method}                & \textbf{Parameter}                    & \textbf{Value}                                 & \textbf{Computation Time (s)} \\ \hline
\multirow{3}{*}{Weighted Sum}  & Weight vectors (\(w_l\))              & \(w_1 = \{0, 0.01, \ldots, 1\}, \ w_2 = 1-w_1\) & \multirow{3}{*}{0.87}         \\ 
                               & Initial guess (\(\mathbf{x}_0\))      & Multiple initial values                       &                              \\ 
                               & Termination criterion (\texttt{xtol}) & \(10^{-8}\)                                   &                              \\ \hline
\multirow{4}{*}{Desirability}  & Lower bound (\(L_j\))                 & \(L_1 = 600 , \ L_2 = 0.6\)                   & \multirow{4}{*}{10.20}        \\ 
                               & Upper bound (\(U_j\))                 & \(U_1 = 725 , \ U_2 = 0.7\)                   &                              \\ 
                               & Shape parameters (\(r_j\))            & \(r_1 = 2.5, r_2 = 0.25\)                     &                              \\
                               & Objective weights (\(\rho_j\))           & \(\rho_1 = \rho_2 = 0.5\)                           &                              \\ \hline
\multirow{2}{*}{NSGA-II}       & Population size (\(N\))               & 300                                           & \multirow{2}{*}{436.73}        \\ 
                               & Generations                          & 1000                                          &                              \\ \hline
\end{tabular}
}
\caption{Hyperparameters, settings, and computation times for the optimization methods for solving Problem I. The computation time for the Weighted Sum Method corresponds to the longest runtime among multiple initial values.}
\label{tab:hyperparameters}
\end{table}

\begin{figure}[htb!]
    \centering
    \includegraphics[width=0.85\textwidth]{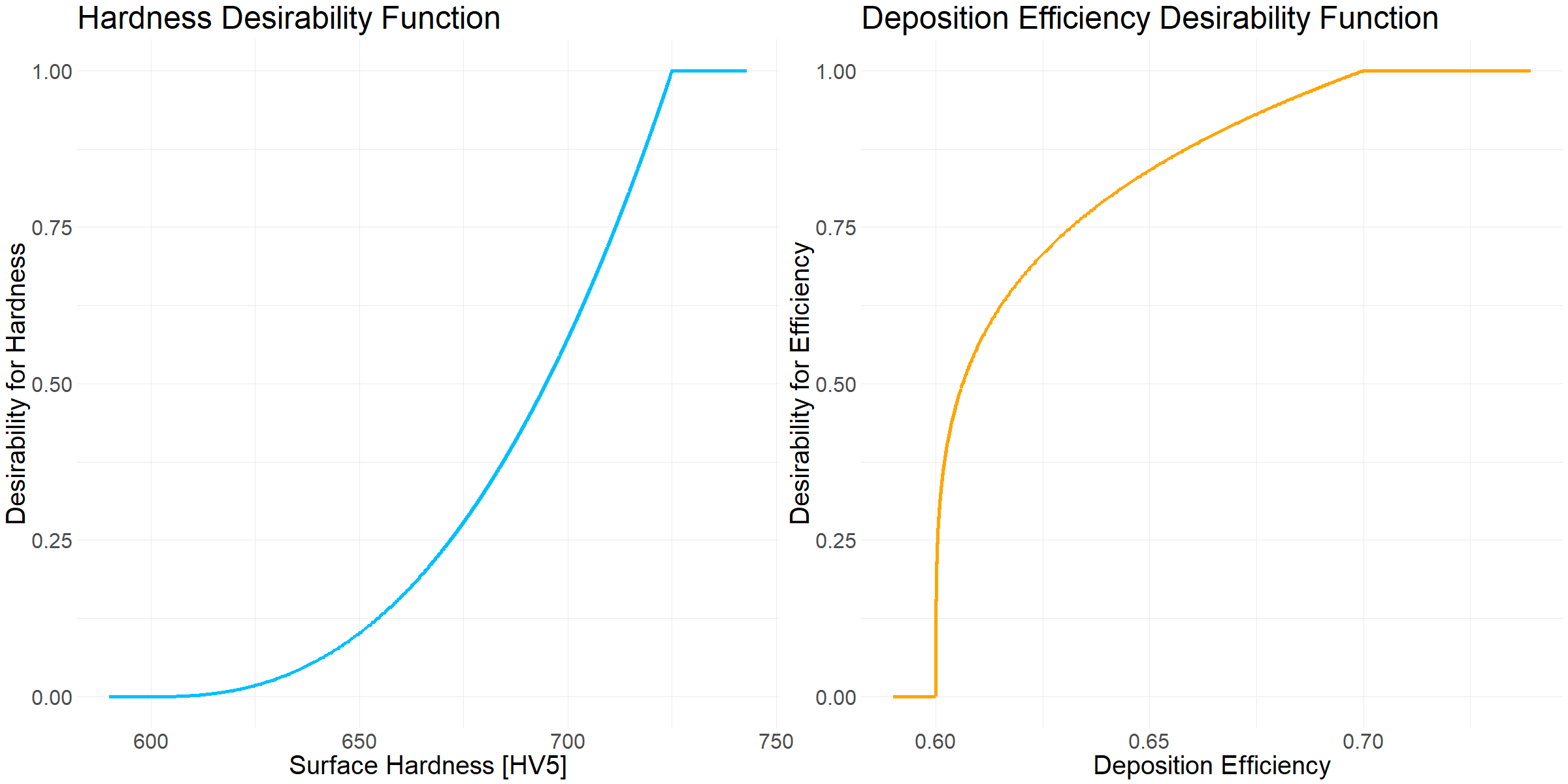}
    \caption{This figure illustrates the one-sided desirability transformations for the dual objectives of maximizing surface hardness (left) and deposition efficiency (right) in the HVOF optimization problem. For surface hardness, the decision maker, based on expert input, defined a lower bound of 600 HV5 and an upper bound of 725 HV5, with a shape parameter $r =2.5$. For deposition efficiency, the lower and upper bounds were set to 0.6 and 0.7, respectively, with a shape parameter $r =0.25$.}
    \label{fig:desirability_plot}
\end{figure}

The results of the dual-objective optimization problem, aimed at maximizing surface hardness and deposition efficiency, provide insights into the trade-offs inherent in the HVOF process. Figure~\ref{fig:desirability_plot} illustrates the desirability functions for both objectives, which form the desirability-based optimization approach. The desirability function for surface hardness, defined with a lower bound of 600 HV5 and an upper bound of 725 HV5, reflects the importance of achieving higher hardness values for industrial applications. Its convex shape penalizes values closer to the lower bound while strongly favoring those near the upper bound. Conversely, the desirability function for deposition efficiency, with bounds of 0.6 and 0.7 and a concave shape, indicates a more lenient optimization approach, emphasizing practicality and economic viability over absolute efficiency.

\begin{figure}[htb!]
    \centering
    \includegraphics[width=0.85\textwidth]{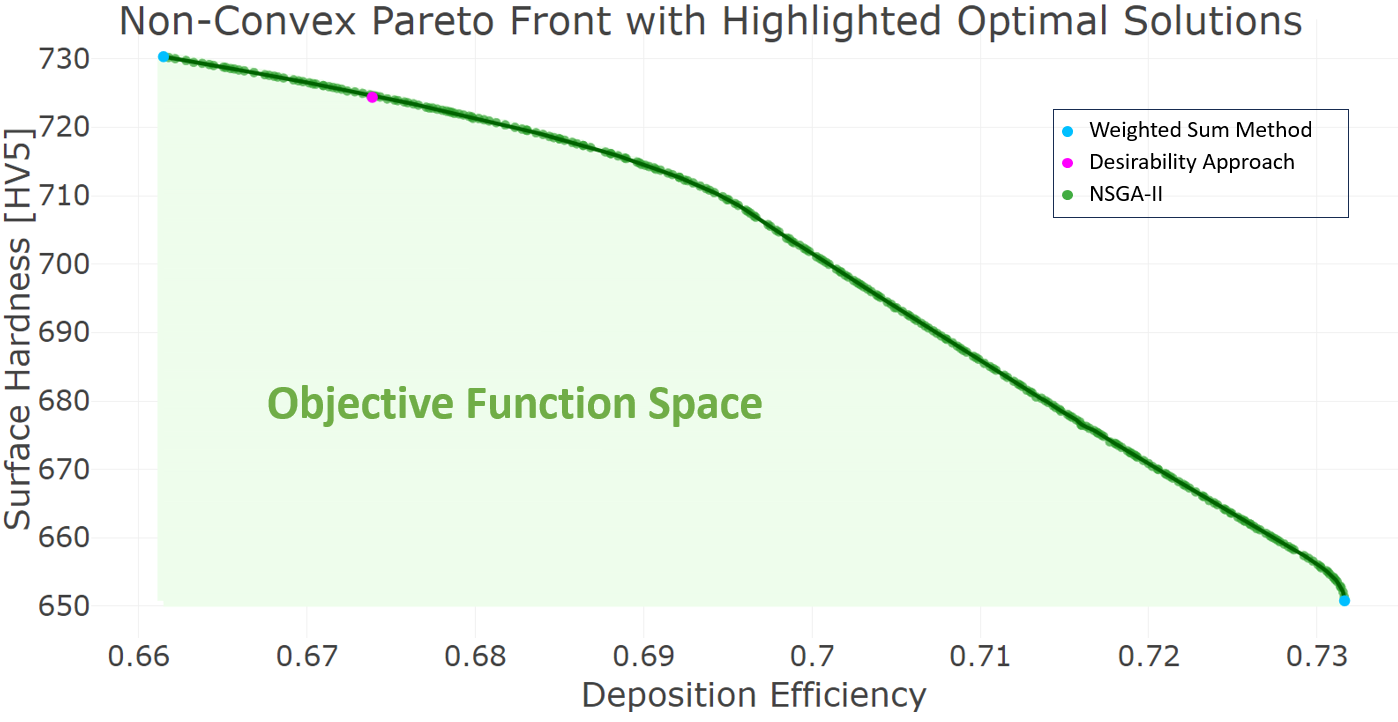}
    \caption{This figure depicts the non-convex Pareto front for the dual-objective optimization problem of maximizing surface hardness and deposition efficiency. The full set of Pareto-optimal solutions, identified by NSGA-II, is represented by green points, with a connecting green line indicating the Pareto front. The magenta point highlights the solution with the highest overall desirability, determined through the desirability approach based on the decision maker's predefined preference criteria. Due to the non-convex nature of the Pareto front, the Weighted Sum Method converges only to extreme solutions, shown as blue points, underscoring its limitations in exploring the entirety of the trade-off surface.}
    \label{fig:non_convexPareto}
\end{figure}

The non-convex Pareto front obtained using NSGA-II is depicted in Figure~\ref{fig:non_convexPareto}, showcasing the full set of Pareto-optimal solutions (green points). These solutions highlight the trade-off between the two objectives, emphasizing the limitations of conventional methods like the weighted sum approach, which identifies solutions only at the extremes of the Pareto front (blue points). The solution with the highest overall desirability, identified based on the decision maker's predefined preference criteria, is marked as a magenta point on the front. This solution serves as a reference point for practical implementation, balancing both objectives according to expert-defined priorities.

The results obtained with NSGA-II emphasize its ability to comprehensively explore the solution space, capturing the trade-offs that gradient-based optimization methods fail to identify. However, the computational expense of NSGA-II compared to the weighted sum or desirability approaches may become a factor for industrial applications where rapid optimization is required. However, in the context of industrial thermal spraying, where the coating of large-scale components typically involves extensive preparation and masking, computational time is a secondary concern relative to overall process efficiency. To provide a complete comparison of the methods, the computation times are summarized in Table~\ref{tab:hyperparameters}. These calculations were conducted on a standard PC equipped with an Intel(R) Core(TM) i7-11850H vPro processor (2.5 GHz, 8 cores), with only a single core used for the serial computations.

Based on the optimization results, the solution identified using the desirability approach is selected for practical implementation due to its alignment with industrial priorities and its balance of the dual objectives. This solution corresponds to a theoretical surface hardness of $724.38$ HV5 and a deposition efficiency of $0.674$. Additionally, an alternative solution from the Pareto front, identified using the NSGA-II algorithm, is selected for further exploration. This solution achieves a theoretical surface hardness of $711.86$ HV5 and a deposition efficiency of $0.693$. The latter solution is chosen to investigate process conditions where deposition efficiency is maximized without significant compromise in surface hardness, as a noticeable drop in hardness occurs when deposition efficiency exceeds $0.69$. A summary of these solutions is presented in Table~\ref{tab:solution_all} in the subsequent Section~\ref{sect_validation}.

\subsection{Problem II: Optimization of Surface Hardness, Deposition Efficiency, and Particle Temperature}

Building on Problem I, this second optimization problem introduces a third objective: minimizing particle temperature. This addition addresses specific industrial requirements, particularly in applications where excessive thermal load on the substrate must be avoided. Excessive particle temperatures can lead to undesired phase transformations within the coating, compromising its performance. Additionally, high thermal loads risk component deformation, especially in temperature-sensitive materials such as aluminum \cite{umantsev2007thermal}.

Incorporating particle temperature as an optimization objective increases the complexity of the problem, requiring a careful balance between coating performance, process efficiency, and thermal constraints. The same multi-objective optimization framework is employed to systematically explore parameter settings that satisfy all three objectives. The optimization problem is formulated as:

\[
\text{Minimize } F(\boldsymbol{x}) = \left[-f_1(\mathbf{x}), -f_2(\mathbf{x}), f_3(\mathbf{x})\right]^\top 
\]
where $f_1(\mathbf{x})$ represents surface hardness, $f_2(\mathbf{x})$ deposition efficiency, and $f_3(\mathbf{x})$ particle in-flight temperature. The negative signs for $f_1$ and $f_2$ indicate maximization, while $f_3$ is minimized directly. This optimization, as before, is subject to the following physical constraints on the process parameters:

\[
x_{\text{Lower}} \leq x \leq x_{\text{Upper}}
\]
where \( x_{\text{Lower}} \) and \( x_{\text{Upper}} \) denote material-specific boundaries on process parameters, as detailed in Table~\ref{tab:limits}.

Table~\ref{tab:hyperparametersII} summarizes the hyperparameters and settings employed for the optimization methods addressing Problem II. In the desirability-based approach, the desirability functions for surface hardness and deposition efficiency remain largely unchanged from Problem I, with a slight adjustment to the deposition efficiency bounds ($L_2=0.5$ and $U_2=0.65$) to reflect refined industrial requirements. The determination of lower and upper bounds for particle temperature ($L_3=1600$ and $U_3=1720$) was informed by domain expertise, prioritizing lower particle temperatures to mitigate substrate deformation and phase transformations. The shape parameters for this desirability function were adjusted to favor even cooler temperatures.

\begin{table}[!ht]
\centering
\resizebox{\textwidth}{!}{
\begin{tabular}{llll}
\hline
\textbf{Method}                & \textbf{Parameter}                    & \textbf{Value}                                 & \textbf{Computation Time (s)} \\ \hline
\multirow{3}{*}{Weighted Sum}  
& Weight vectors (\(w_l\))     & \((w_1, w_2, w_3) \) & \multirow{3}{*}{11.30} \\ 
                               & Initial guess (\(\mathbf{x}_0\))      & Multiple initial values                           &                              \\ 
                               & Termination criterion (\texttt{xtol}) & \(10^{-8}\)                                   &                              \\ \hline
\multirow{4}{*}{Desirability}  & Lower bound (\(L_j\))                 & \(L_1 = 600 , \ L_2 = 0.5, \ L_3 = 1600\)                   & \multirow{4}{*}{11.97}        \\ 
                               & Upper bound (\(U_j\))                 & \(U_1 = 725 , \ U_2 = 0.65, \ U_3 = 1720\)                   &                              \\ 
                               & Shape parameters (\(r_j\))            & \(r_1 = 2.5, r_2 = 0.25, r_3 = 2\)                     &                              \\
                               & Objective weights (\(\rho_j\))           & \(\rho_1 = \rho_2 = \rho_3 = 1/3\)                           &                              \\ \hline
\multirow{2}{*}{NSGA-II}       & Population size (\(N\))               & 5000                                           & \multirow{2}{*}{9262.75}        \\ 
                               & Generations                          & 100                                          &                              \\ \hline
\end{tabular}
}
\caption{Hyperparameters, settings, and computation times for the optimization methods for solving Problem II. The weight vector \((w_1, w_2, w_3)\) in the Weighted Sum Method is defined such that \(w_1, w_2, w_3 \geq 0\), \(\sum_{j=1}^3 w_j = 1\), \({w_1, w_2 \in \{0, 0.01, \ldots, 1\}}\), and \(w_3 = 1 - (w_1 + w_2)\). The computation time for the Weighted Sum Method corresponds to the longest runtime among multiple initial values.}
\label{tab:hyperparametersII}
\end{table}

For the Weighted Sum Method, an additional weight ($w_3$) was introduced for particle temperature. The weight vectors were systematically designed to satisfy the imposed requirements, i.e., 
\begin{align*}
    (w_1, w_2, w_3 \geq 0), \ \sum_{j=1}^3 w_j = 1, \quad w_1, w_2 \in \{0, 0.01, \ldots, 1\}, \ w_3 = 1 - (w_1 + w_2). 
\end{align*}

For the NSGA-II method, the population size was increased to $N=5000$ to enhance the representation of the Pareto front, while the number of generations was reduced to 100 for computational efficiency. The overall objective remains the simultaneous maximization of deposition efficiency and surface hardness while minimizing particle temperature. 

\begin{figure}[htb!]
    \centering
    \includegraphics[width=0.85\textwidth]{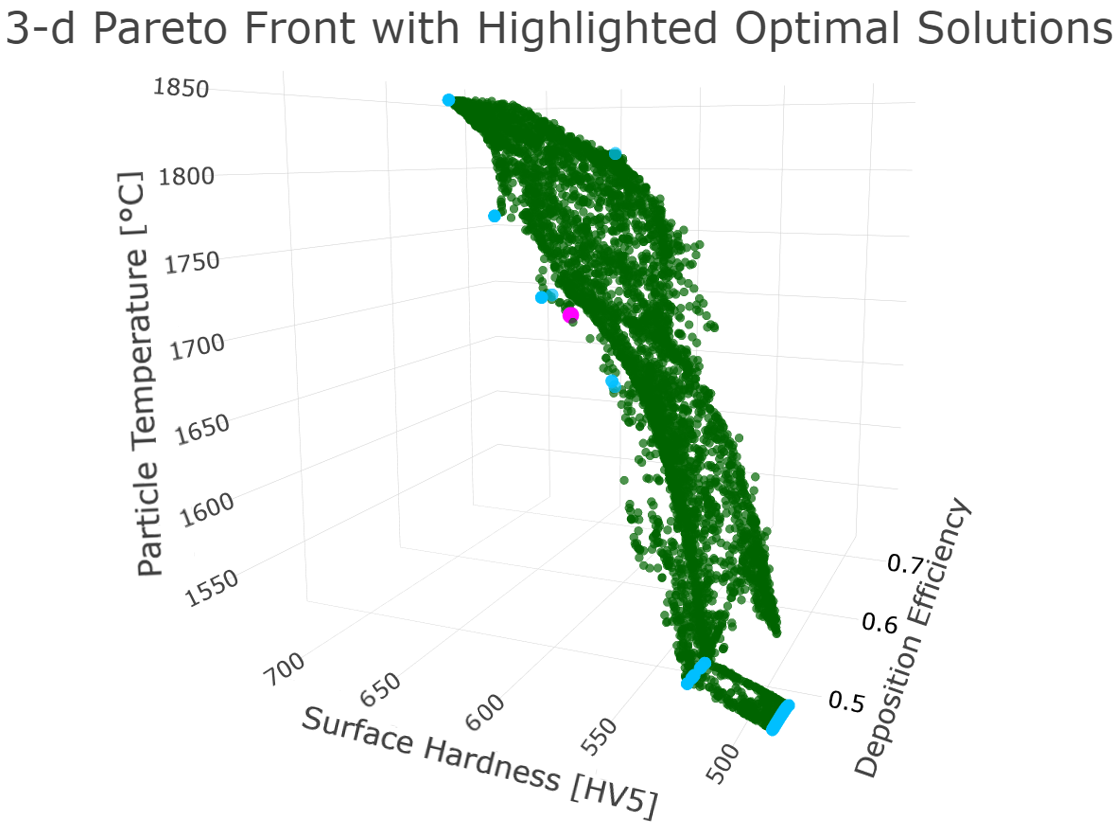}
    \caption{Three-dimensional Pareto front for Problem II, illustrating the trade-offs among deposition efficiency, surface hardness (HV5), and particle temperature (°C). The Pareto-optimal solutions identified by NSGA-II are represented as green points. The magenta point highlights the solution with the highest overall desirability, determined using the desirability approach based on the decision maker's predefined preference criteria. In contrast, the Weighted Sum Method (blue points) converges only to extreme solutions at the boundaries of the Pareto front.}
    \label{fig:II_pareto}
\end{figure}

Figure~\ref{fig:II_pareto} presents the three-dimensional Pareto front for Problem II, showcasing the trade-offs among deposition efficiency, surface hardness, and particle temperature. The set of Pareto-optimal solutions identified by the NSGA-II algorithm is depicted as green points, mapping the diverse trade-offs between the three objectives. The magenta point represents the solution with the highest overall desirability, determined using the desirability approach. This solution aligns with the decision maker's predefined preference criteria, emphasizing a balanced compromise among the competing objectives.

In contrast, the solutions identified by the Weighted Sum Method, represented as blue points, are concentrated at the boundaries of the Pareto front. This behavior reflects the limitation of the method in exploring intermediate trade-offs when objectives are non-convex \cite{das1997closer,nakayama1995aspiration}. Notably, since no additional constraints were imposed during optimization, both NSGA-II and the Weighted Sum Method also identify solutions which fall outside the range of preferred outcomes defined by the desirability approach. The advantage of employing an a posteriori approach, as exemplified by NSGA-II, lies in its ability to provide a comprehensive set of trade-off solutions, enabling the selection of a preferred solution post-optimization. The computation times for the methods are detailed in Table~\ref{tab:hyperparametersII}, highlighting the differences between the exhaustive exploration of the Pareto front (NSGA-II) and the relatively rapid convergence of the desirability and weighted sum approaches.

Building on the optimization results, two solutions were identified by application domain experts as the most promising candidates for further evaluation. These solutions, summarized in Table~\ref{tab:solution_all}, represent a selection of the optimal trade-offs among the three objectives and are evaluated through practical validation in Section~\ref{sect_validation} to assess the accuracy and reliability of the GLM models in predicting real-world performance.

\subsection{Problem III: Optimization of Coating Porosity, Coating Roughness, and Particle Temperature}

The third optimization problem addresses the simultaneous optimization of coating properties which are not specific to tungsten carbide cobalt chrome but are critical across various industrial applications. The objective is to achieve a highly dense coating (minimizing porosity) with an exceptionally smooth surface (minimizing roughness) while limiting thermal stress during the process (minimizing particle in-flight temperature). While post-processing techniques can refine surface roughness, achieving a smoother as-sprayed coating provides an advantageous starting point.

In industrial applications such as wear-resistant rolls in steel manufacturing, these properties are particularly relevant. Low porosity enhances durability and resistance to mechanical stress-induced cracking, while reduced roughness minimizes friction during rolling operations. Additionally, maintaining low particle temperatures is essential to prevent thermal distortion or metallurgical changes which could compromise roll performance and lifespan \cite{umantsev2007thermal}. The optimization problem is formulated as:

\[
\text{Minimize } F(\boldsymbol{x}) = \left[f_1(\mathbf{x}), f_2(\mathbf{x}), f_3(\mathbf{x})\right]^\top 
\]
where $f_1(\mathbf{x})$, $f_2(\mathbf{x})$, $f_3(\mathbf{x})$ correspond to coating porosity, surface roughness, and particle in-flight temperature, respectively. The optimization is constrained by:
\[
x_{\text{Lower}} \leq x \leq x_{\text{Upper}}
\]
where \( x_{\text{Lower}} \) and \( x_{\text{Upper}} \) denote the material-specific bounds for the process input parameters, as specified in Table~\ref{tab:limits}.

The optimization framework remains consistent with Problem II, with modifications limited to objective-specific desirability limits and shape parameters to accommodate the distinct optimization targets. The desirability bounds ($L_j$,$U_j$) and shape parameters ($r_j$) were adjusted to reflect the target ranges for porosity, roughness, and particle temperature, as detailed in Table~\ref{tab:hyperparametersIII}.

The weight vectors \(w_l\) for the Weighted Sum Method were constructed following the same systematic approach as in Problem II. Furthermore, the optimization is initialized using multiple starting values $\boldsymbol{x}_0$ for the decision variables with a termination criterion of \texttt{xtol}\(=10^{-8}\). The NSGA-II settings, including a population size of $N$=5000 and 100 generations, were retained from Problem II to ensure methodological consistency. Table~\ref{tab:hyperparametersIII} provides a summary of the hyperparameters, settings, and computational effort for each optimization method applied to Problem III.

\begin{table}[ht]
\centering
\resizebox{\textwidth}{!}{
\begin{tabular}{llll}
\hline
\textbf{Method}                & \textbf{Parameter}                    & \textbf{Value}                                 & \textbf{Computation Time (s)} \\ \hline
\multirow{3}{*}{Weighted Sum}  
& Weight vectors (\(w_l\))     & \((w_1, w_2, w_3) \) & \multirow{3}{*}{11.15} \\ 
                               & Initial guess (\(\mathbf{x}_0\))      & Multiple initial values                           &                              \\ 
                               & Termination criterion (\texttt{xtol}) & \(10^{-8}\)                                   &                              \\ \hline
\multirow{4}{*}{Desirability}  & Lower bound (\(L_j\))                 & \(L_1 = 13 , \ L_2 = 26, \ L_3 = 1600\)                   & \multirow{4}{*}{12.69}        \\ 
                               & Upper bound (\(U_j\))                 & \(U_1 = 15 , \ U_2 = 35, \ U_3 = 1720\)                   &                              \\ 
                               & Shape parameters (\(r_j\))            & \(r_1 = 1.5, r_2 = 2.5, r_3 = 2\)                     &                              \\
                               & Objective weights (\(\rho_j\))           & \(\rho_1 = \rho_2 = \rho_3 = 1/3\)                           &                              \\ \hline
\multirow{2}{*}{NSGA-II}       & Population size (\(N\))               & 5000                                           & \multirow{2}{*}{10026.14}        \\ 
                               & Generations                          & 100                                          &                              \\ \hline
\end{tabular}
}
\caption{Hyperparameters, settings, and computation times for the optimization methods for solving Problem III. The weight vector \((w_1, w_2, w_3)\) in the Weighted Sum Method is defined such that \(w_1, w_2, w_3 \geq 0\), \(\sum_{j=1}^3 w_j = 1\), \({w_1, w_2 \in \{0, 0.01, \ldots, 1\}}\), and \(w_3 = 1 - (w_1 + w_2)\). The computation time for the Weighted Sum Method corresponds to the longest runtime among multiple initial values.}
\label{tab:hyperparametersIII}
\end{table}

\begin{figure}[htb!]
    \centering
    \includegraphics[width=1.00\textwidth]{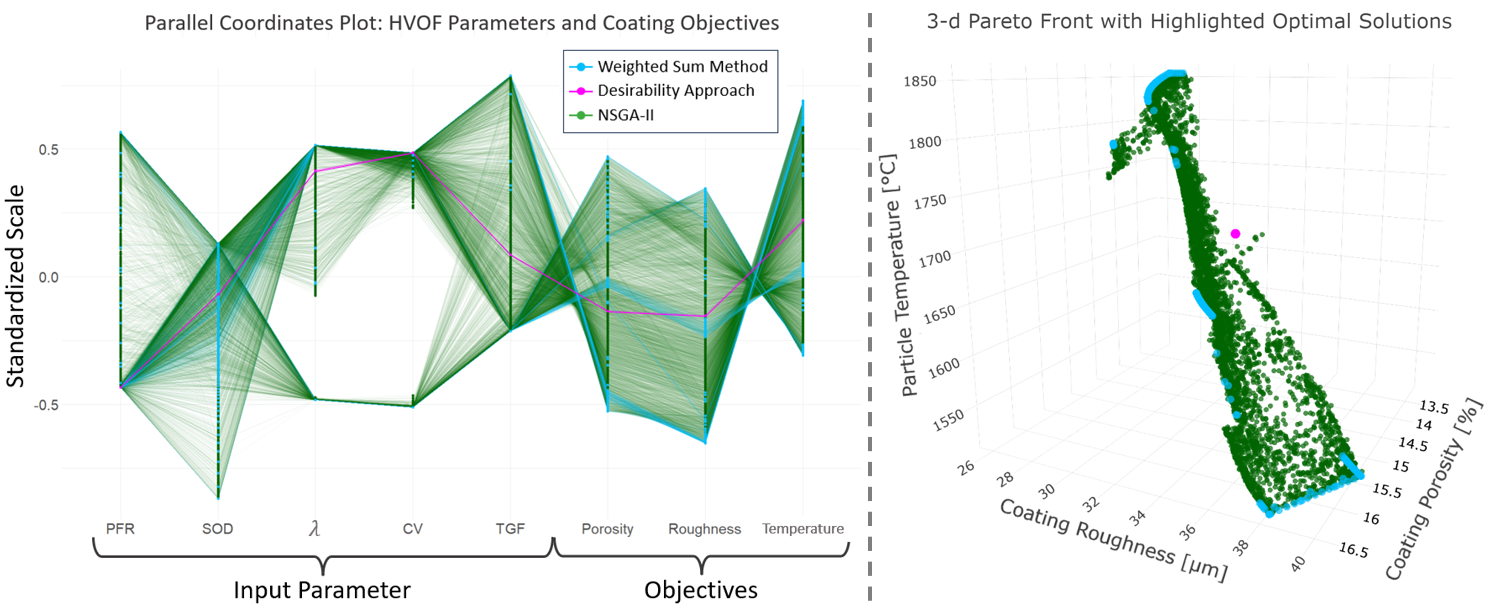}
    \caption{Left: Alternative illustration of the optimization results for Problem III using a parallel coordinates plot. The first five vertical axes represent the normalized HVOF input parameters (powder feed rate (PFR), stand-off distance (SOD), fuel-to-oxygen ratio ($\lambda$), coating velocity (CV), and total gas flow (TGF)), while the last three axes correspond to the scaled optimization objectives (porosity, roughness, and particle temperature). Each path represents a solution, with blue lines indicating results from the Weighted Sum Method with SQP, the magenta line highlighting the solution with the highest desirability, and green lines depicting Pareto-optimal solutions identified by NSGA-II. Right: Three-dimensional Pareto front for Problem III, illustrating the trade-offs among coating porosity (\%), coating roughness (µm), and particle temperature (°C).}
    \label{fig:III_parallelPlot}
\end{figure}

The results of Problem III are visualized in Figure~\ref{fig:III_parallelPlot} using a parallel coordinates plot (left). This alternative representation highlights the relationships between the normalized HVOF input parameters (powder feed rate, stand-off distance, fuel-to-oxygen ratio, coating velocity, and total gas flow) and the scaled optimization objectives (porosity, roughness, and particle temperature). The solutions are represented as paths, with blue lines corresponding to the Weighted Sum Method using SQP, the magenta line identifying the solution with the highest desirability, and green lines indicating the Pareto-optimal solutions obtained via NSGA-II.

Distinct clusters are evident in the fuel-to-oxygen ratio and coating velocity axes, indicating two distinct parameter configurations for achieving optimal outcomes. Specifically, lower fuel-to-oxygen ratios are paired with lower coating velocities, while higher fuel-to-oxygen ratios are coupled with higher coating velocities. To achieve a dense coating, higher total gas flow values are generally required, which aligns with expectations. The Weighted Sum Method identified several solutions, with many located at the boundaries of the objectives, indicating extreme trade-offs. However, in contrast to previous optimization problems, the Weighted Sum Method in this case shows a broader range of solutions, as evidenced by the blue paths which extend across the interior regions of the last three vertical axes in the parallel coordinates plot. Despite this apparent diversity, it is important to note that these solutions still correspond to the boundaries of the objective space when considered in the context of the optimization problem, as illustrated in Figure~\ref{fig:III_parallelPlot} (right).

From the optimization results for Problem III, two solutions were chosen by experts in the application domain as the most suitable candidates for further investigation. These solutions, detailed in Table~\ref{tab:solution_all}, illustrate a selection of optimal compromises among the three objectives: coating porosity, coating roughness, and particle temperature. Their practical feasibility will be evaluated in the subsequent Section~\ref{sect_validation}, where the predictive capabilities of the GLMs will be examined, and coatings will be performed using the corresponding HVOF input parameters to validate the optimization outcomes.

\section{Validation of Optimization Solutions}\label{sect_validation}
To assess the reliability and practical relevance of the optimization results obtained in Section~\ref{sect_application}, validation trials were conducted. In particular, the objective was to verify whether the optimized parameter settings yield the expected coating properties under real operating conditions. Specifically, the input parameter configurations corresponding to the two selected solutions for each optimization problem were implemented in the spray booth. Samples produced under these conditions were analyzed to compare the predicted properties with the actual outcomes. This evaluation aims to determine the degree of agreement between the theoretical solutions and their practical realizations, providing insight into the robustness and applicability of the optimization approach. 

The validation trials were conducted using the coating material (WC-10Co-4Cr) and the HVOF coating system detailed in Section~\ref{sect_background}. Table~\ref{tab:solution_all} summarizes the optimization results, presenting the input parameter configurations selected for validation by expert technicians for the three optimization problems introduced in Section~\ref{sect_application}. These configurations specify the process settings for powder feed rate (PFR), stand-off distance (SOD), fuel-to-oxygen ratio ($\lambda$), coating velocity (CV), and total gas flow (TGF), alongside the corresponding predicted values for the objectives. Recall that the objectives vary across problems, with Problem I targeting maximum surface hardness and deposition efficiency, Problem II introducing particle temperature as a third objective to minimize, and Problem III emphasizing the minimization of coating porosity, surface roughness, and particle temperature.

\begin{table}[ht]
\footnotesize
\centering
\renewcommand{\arraystretch}{1.2}
\resizebox{\textwidth}{!}{
\begin{tabular}{l|l|ccccc|ccc}
\hline\hline
\textbf{Problem} & \textbf{Solution} & \textbf{PFR} & \textbf{SOD} & \boldmath{$\lambda$} & \textbf{CV} & \textbf{TGF} & \textbf{$f_1(\mathbf{x})$} & \textbf{$f_2(\mathbf{x})$} & \textbf{$f_3(\mathbf{x})$} \\ \hline
\multirow{2}{*}{\textbf{I}} 
& Desirability & \(45.00\) & \(200.00\) & \(1.04\) & \(80.00\) & \(751.00\) & \(724.38 \ HV5\) & \(67.4 \ \%\) & -- \\ 
& NSGA-II & \(45.00\) & \(200.00\) & \(1.04\) & \(90.67\) & \(751.00\) & \(711.86 \ HV5\) & \(69.3 \ \%\) & -- \\ \hline

\multirow{2}{*}{\textbf{II}} 
& NSGA-II I & \(49.10\) & \(259.22\) & \(0.84\) & \(101.01\) & \(727.73\) & \(604.71 \ HV5\) & \(59.5 \ \%\) & \(1690.97\) °C \\ 
& NSGA-II II & \(62.00\) & \(259.74\) & \(0.85\) & \(99.71\) & \(748.74\) & \(603.73 \ HV5\) & \(61.7 \ \%\) & \(1698.52\) °C \\ \hline

\multirow{2}{*}{\textbf{III}} 
& NSGA-II I & \(45.02\) & \(259.96\) & \(1.04\) & \(120.01\) & \(638.88\) & \(14.57\ \%\) & \(32.71\) µm & \(1664.19\) °C\\ 
& NSGA-II II & \(45.00\) & \(255.53\) & \(1.04\) & \(118.61\) & \(615.34\) & \(14.95\ \%\) & \(34.10\) µm & \(1622.64\) °C\\ \hline \hline
\end{tabular}
}
\caption{Optimization solutions for three problems, detailing input parameters and corresponding objectives. Problem I targets maximizing $f_1(\mathbf{x}) =$ surface hardness and $f_2(\mathbf{x}) =$ deposition efficiency, Problem II consists of minimizing $f_3(\mathbf{x}) =$ particle temperature as a third objective, and Problem III focuses on minimizing $f_1(\mathbf{x}) =$ coating porosity, $f_2(\mathbf{x})=$ roughness, and $f_3(\mathbf{x}) =$ particle temperature. Abbreviations: PFR = Powder Feed Rate, SOD = Stand-Off Distance, $\lambda$ = Fuel-to-Oxygen Ratio, CV = Coating Velocity, TGF = Total Gas Flow.}
\label{tab:solution_all}
\end{table}

\begin{figure}[htb!]
    \centering
    \includegraphics[width=0.7\textwidth]{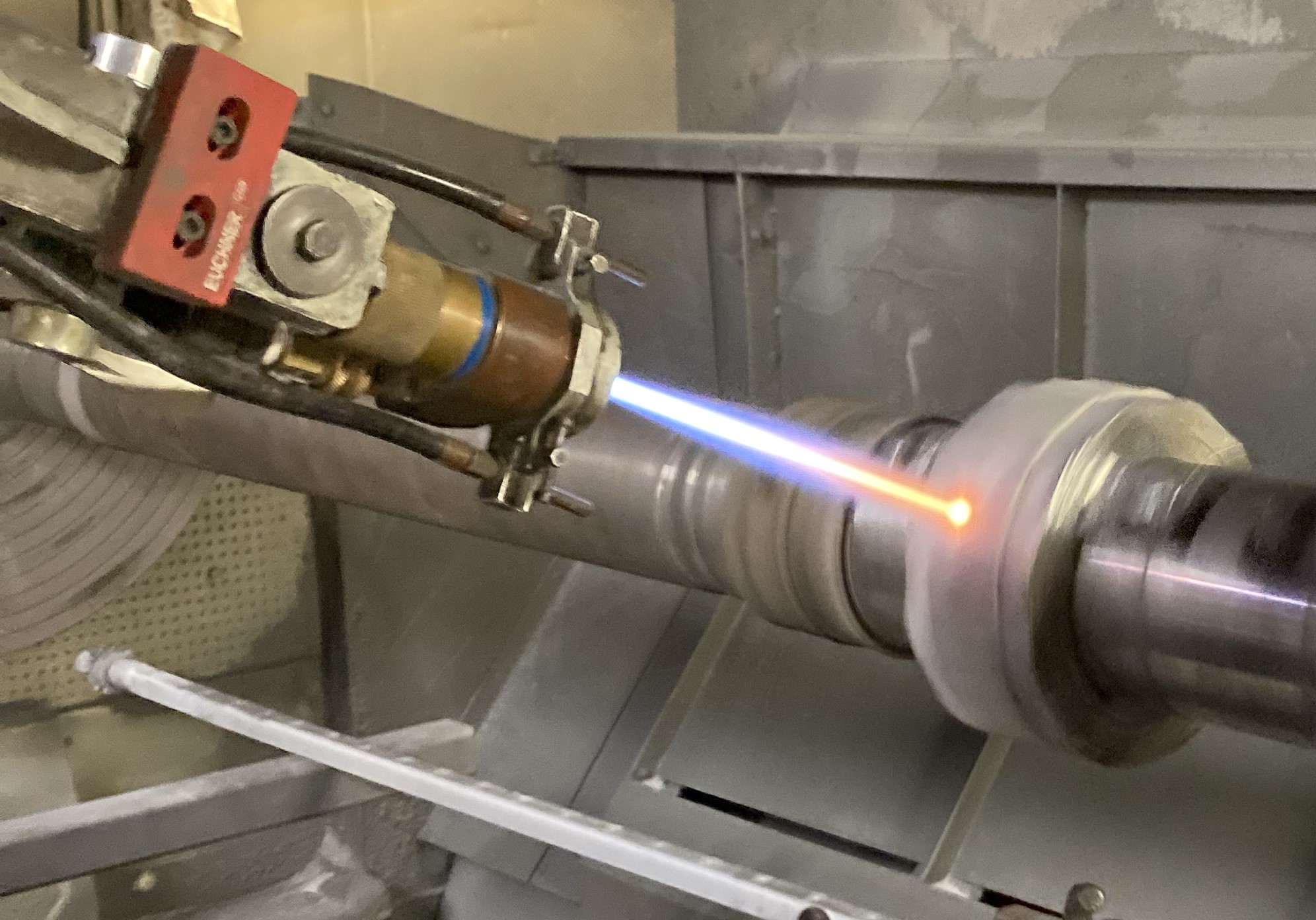}
    \caption{Photograph of the experimental setup used for the practical validation of the optimization results in the HVOF coating process. The image highlights key components, including the HVOF torch, turning lathe, sample holder, and the active coating stream.}
    \label{fig:HVOF_real}
\end{figure}

\begin{figure}[htb!]
    \centering
    \includegraphics[width=0.85\textwidth]{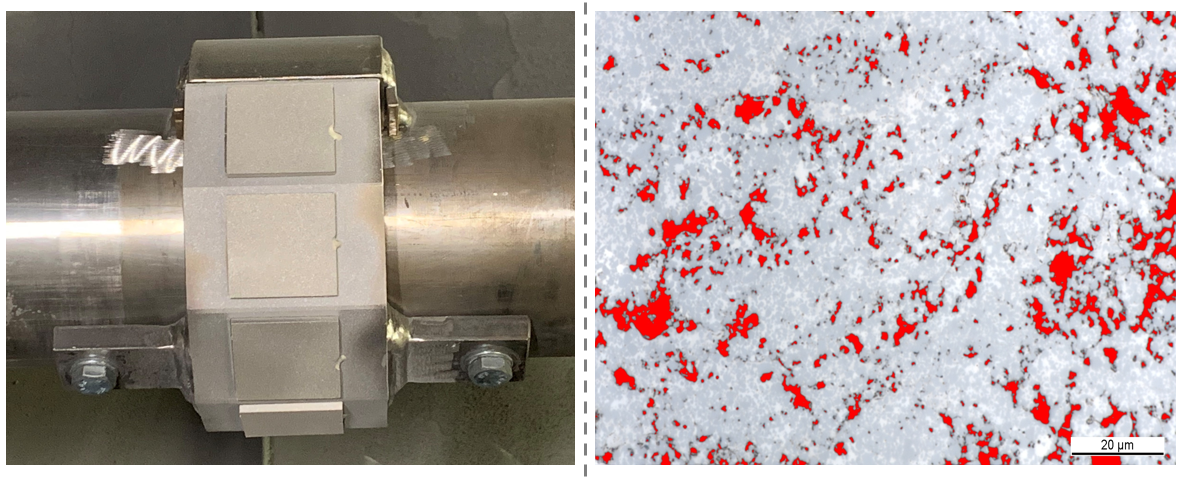}
    \caption{Left: Photograph of the sample holder with mounted and pretreated sample plates, prepared for the HVOF coating process after activation via grit blasting. Right: Microscopic image illustrating the porosity quantification technique using image analysis \cite{ang2014review}. The method highlights void areas in red, enabling precise calculation of the porosity percentage for the specific sample.}
    \label{fig:sapleholder_porosity}
\end{figure}

The experimental setup used for the validation trials is illustrated in Figure~\ref{fig:HVOF_real}, showcasing key components such as the HVOF torch, turning lathe, sample holder, and the active coating stream. The sample plates, prepared and mounted on the sample holder, underwent grit blasting to ensure proper surface activation prior to coating (Figure~\ref{fig:sapleholder_porosity}, left). This surface preparation step is critical for achieving optimal adhesion and coating quality.

Quantification of porosity was performed using an image analysis technique \cite{ang2014review}, as shown in Figure~\ref{fig:sapleholder_porosity} (right). Microscopic imaging of the coated samples identifies void regions, which are highlighted in red, enabling precise calculation of the porosity percentage by determining the proportion of voids relative to the total sample area. The final porosity value was obtained by averaging five measurements taken from different regions of each sample. Particle temperature was measured in-flight using a Spraywatch camera with the SW4 software suite (Oseir), ensuring accurate monitoring of thermal conditions during the coating process. The reported temperature represents the average value over the entire coating duration. Surface hardness was evaluated on the coated samples using a Cisam-Ernst S.r.l E-Computest mobile hardness tester, employing a spheroconical diamond indenter at a load of 5 kg and a testing time of 2 seconds. Following the approach outlined in \cite{mathesius2009praxis}, ten hardness measurements were taken per sample, with the highest and lowest values discarded before averaging the remaining data. Deposition efficiency was determined by weighing the sample before and after the coating process and correcting for the total material input. Together, these measurement techniques provide a comprehensive assessment of the coating properties, facilitating the validation of the optimization results.

The validation trials provide a quantitative assessment of the optimization solutions derived in Section~\ref{sect_application}, enabling a direct comparison between theoretical predictions and experimentally observed coating properties. Table~\ref{tab:validation_results} summarizes the deviations between the predicted and measured values for each optimization problem, offering insights into the accuracy and practical feasibility of the applied optimization methods. 

\begin{table}[ht]
\centering
\footnotesize
\renewcommand{\arraystretch}{1.2}
\resizebox{\textwidth}{!}{
\begin{tabular}{l|l|l|c|c|c}
\hline\hline
\textbf{Problem} & \textbf{Solution} & \textbf{Property} & \textbf{Theoretical Value} & \textbf{Practical Value} & \textbf{Deviation (\%)} \\ \hline
\multirow{4}{*}{\textbf{I}} 
& Desirability & Surface Hardness (HV5) & \(724.38\) & \(719.50\) & \(-0.67\) \\ 
&              & Deposition Efficiency (\%) & \(67.4\) & \(68.7\) & \(+1.93\) \\ 
& NSGA-II      & Surface Hardness (HV5) & \(711.86\) & \(709.00\) & \(-0.40\) \\ 
&              & Deposition Efficiency (\%) & \(69.3\) & \(69.4\) & \(+0.14\) \\ \hline

\multirow{6}{*}{\textbf{II}} 
& NSGA-II I    & Surface Hardness (HV5) & \(604.71\) & \(608.10\) & \(+0.56\) \\ 
&              & Deposition Efficiency (\%) & \(59.5\) & \(60.1\) & \(+1.01\) \\ 
&              & Particle Temperature (°C) & \(1690.97\) & \(1694.11\) & \(+0.19\) \\ 
& NSGA-II II    & Surface Hardness (HV5) & \(603.73\) & \(601.00\) & \(-0.45\) \\ 
&              & Deposition Efficiency (\%) & \(61.7\) & \(62.4\) & \(+1.13\) \\ 
&              & Particle Temperature (°C) & \(1698.52\) & \(1696.14\) & \(-0.14\) \\ \hline

\multirow{6}{*}{\textbf{III}} 
& NSGA-II I    & Porosity (\%) & \(14.57\) & \(14.52\) & \(-0.34\) \\ 
&              & Roughness (\(\mu\)m) & \(32.71\) & \(32.19\) & \(-1.59\) \\ 
&              & Particle Temperature (°C) & \(1664.19\) & \(1672.75\) & \(+0.51\) \\ 
& NSGA-II II   & Porosity (\%) & \(14.95\) & \(13.74\) & \(-8.09\) \\ 
&              & Roughness (\(\mu\)m) & \(34.10\) & \(34.94\) & \(+2.46\) \\ 
&              & Particle Temperature (°C) & \(1622.64\) & \(1610.37\) & \(-0.76\) \\ \hline\hline
\end{tabular}
}
\caption{Comparison of theoretical and practical results for optimization solutions across three problems. Deviations (\%) are calculated as \((\text{Practical} - \text{Theoretical}) / \text{Theoretical} \times 100\).}
\label{tab:validation_results}
\end{table}

The validation results indicate a high degree of agreement between the theoretical predictions and practical outcomes, with deviations remaining within a small percentage range across all optimization problems. For Problem I, the surface hardness predictions show a deviation of less than 1\%, with the NSGA-II solution achieving an error margin of only 0.40\%. Similarly, deposition efficiency reveals minimal discrepancies, with a maximum deviation of 1.93\%. These findings confirm that the applied models accurately capture the relationship between HVOF process parameters and coating properties, highlighting the suitability of the gamma regression-based formulations for predictive optimization.

For Problem II, where particle temperature is introduced as an additional objective, the applied optimization methodology still demonstrate strong performance. The deviations for surface hardness and deposition efficiency remain below 1.13\%, while particle temperature exhibits an error below 0.2\%. These results suggest that the optimization approach effectively accounts for the third objective without substantially compromising the accuracy of the first two objectives. The slightly increased deviations compared to Problem I may be attributed to the additional complexity introduced by the competing trade-offs among the three objectives.

For Problem III, the measured values for roughness and particle temperature align with the theoretical solutions, with deviations remaining below 2.5\%. However, a more pronounced discrepancy is observed in the porosity measurements, particularly for the NSGA-II II solution, which deviates by approximately 8.09\%. This deviation is notably larger than those observed for the other optimization objectives and can likely be attributed to the known limitations of the image-based porosity measurement technique \cite{ang2014review}. As discussed in \cite{ang2014review}, the accuracy of this method is influenced by several factors, including metallographic sample preparation, imaging technique, and post-processing procedures such as thresholding. These methodological sensitivities introduce variability in the measured porosity values, potentially leading to systematic under- or overestimation. 

The overall accuracy observed across all optimization problems emphasizes the reliability of the proposed methodology. The minor deviations can largely be attributed to the well-adjusted gamma regression models \cite{rannetbauer2024predictive}, which correctly capture nonlinear dependencies within the process. Nevertheless, slight discrepancies may result from experimental uncertainties, including variations in substrate surface conditions, environmental influences, and minor fluctuations in process stability. Despite these factors, the results confirm that the optimization framework provides solutions which are both theoretically valid and practically implementable. This demonstrates its suitability for optimizing coating characteristics while accounting for the practical limitations of the coating system.

\section{Conclusion}\label{sect_conclusion}
This study presents a systematic framework for optimizing critical target variables in HVOF thermal spray processes. By leveraging multi-objective optimization techniques in conjunction with gamma regression models as objective functions, optimal parameter configurations for industrial coating applications were determined. The theoretical foundations and practical implementation of multi-objective optimization were discussed, and a comparative analysis of three distinct optimization approaches highlighted their respective strengths and limitations in HVOF coating optimization.

The desirability-based optimization approach proved intuitive and computationally efficient, making it well-suited for low-dimensional problems. However, its dependence on predefined desirability functions demands precise goal specification, requiring process expertise. The Weighted Sum Method, applied in an a posteriori manner, approximated the Pareto front but exhibited well-documented limitations in handling non-convex objective spaces. Despite multiple initialization strategies, only extreme solutions at the Pareto front boundaries were identified. In contrast, NSGA-II demonstrated the most robust performance, effectively capturing a diverse set of Pareto-optimal solutions. Its primary drawback—computational intensity—is of limited concern in industrial thermal spraying, where preparation and application times outweigh optimization runtimes.

Beyond enforcing box constraints, ensuring the optimization remained within practical and technical system limits was a key consideration. This study extends prior research on predictive modelling of HVOF coatings \cite{rannetbauer2024predictive}, demonstrating the feasibility of GLM-based optimization. While the findings pertain to a specific WC-10Co-4Cr powder and industrial HVOF system, the methodology is generalizable to other coating materials. However, the results underscore the need for domain expertise, as data-driven optimization alone is insufficient without process-specific knowledge.

Future work will focus on extending the proposed framework to alternative coating materials and their associated properties. A key objective is to refine optimization strategies by integrating coating process data with production data from steel manufacturing, thereby enabling a more application-specific and performance-oriented approach. The evolving demands for higher-performance steels, driven by sustainability efforts and stricter CO\textsubscript{2} emission regulations, require coatings capable of resisting new types of operational stresses \cite{rannetbauer2024leveraging}. This shift introduces new challenges in surface engineering, requiring tailored optimization strategies to ensure enhanced durability and functionality.

Furthermore, ongoing research focuses on extending the framework through transfer learning, enabling knowledge transfer from small-scale experimental trials to large-scale industrial applications. Given the lack of reliable non-destructive testing methods for certain coating properties, inductive transfer learning approaches will also be explored to infer unmeasured characteristics from related target variables.

Overall, this study provides a data-driven, multi-objective optimization framework tailored for HVOF thermal spray coatings. By systematically integrating statistical modelling and optimization techniques, it enhances process understanding and facilitates improved control over coating properties. The findings not only advance surface engineering methodologies but also provide a structured framework for industry practitioners and researchers to optimize coating performance under realistic operational constraints, contributing to the broader field of advanced manufacturing.

\section*{Declarations}

\subsection*{Acknowledgements}
    The authors thank voestalpine Stahl GmbH for their support, including the provision of materials, access to their research center, and financial assistance for this study.

\subsection*{Funding}
    This research was funded in part by the Austrian Science Fund (FWF) SFB 10.55776/F68 ``Tomography Across the Scales'', project F6805-N36 (Tomography in Astronomy). For open access purposes, the author has applied a CC BY public copyright license to any author-accepted manuscript version arising from this submission.

\subsection*{Availability of data and materials}
    The datasets generated and analyzed in this study are not publicly accessible due to confidentiality agreements with the company. However, the corresponding author can provide access to the data upon reasonable request. Access will be granted at the discretion of the company and in accordance with confidentiality agreements, on a case-by-case basis.

\subsection*{Competing interests}
     The authors declare that they have no competing interests.

\subsection*{Authors' contributions}
    WR conceived and designed the study, collected the data, performed the analysis, developed the models, applied the optimization algorithms, and wrote the manuscript. WR and CH conducted the experiments, with CH also assisting in data interpretation, providing technical support, and contributing critical revisions. SH and RR offered theoretical insights into the optimization algorithms, recommended methods for detailed analysis, and provided additional critical revisions. All authors reviewed and approved the final manuscript.

\clearpage

\printglossary[type=\acronymtype]
\clearpage
\bibliographystyle{plain}
{\footnotesize
\bibliography{mybib}
}

\newpage
\appendix
\section{Explicit Gamma Regression Model Formulas} \label{appendixA}
In this section, the explicit models for eight key HVOF coating properties are presented, derived from \cite{rannetbauer2024predictive}. Each model expresses the respective coating property as a function of the five input parameters: powder feed rate (PFR), stand-off distance (SOD), fuel-to-oxygen ratio ($\lambda$), coating velocity (CV), and total gas flow (TGF). These models serve as the basis for the application-specific optimization conducted in this work.

\begin{enumerate}
    \item Particle in-flight properties
    \begin{itemize}
        \item Particle in-flight velocity \footnotesize
        \vspace{-0.2cm}
            \begin{align}
            \begin{split}
                f(\text{velocity} | \text{PFR}, \text{SOD}, \lambda, \text{CV}, \text{TGF}) &=     \exp\Big( \beta_0 
                    + \beta_1 \cdot \text{PFR} 
                    + \beta_2 \cdot \text{SOD} 
                    + \beta_3 \cdot \lambda 
                    + \beta_4 \cdot \text{TGF} \\
                    &+ \beta_5 \cdot \text{SOD}^2 
                    + \beta_{6} \cdot \text{TGF}^2 
                    + \beta_{7} \cdot \text{PFR} \cdot \text{SOD} \\
                    &+ \beta_{8} \cdot \text{PFR} \cdot {\lambda}
                    + \beta_{9} \cdot \text{SOD} \cdot \text{TGF} \Big)
            \end{split}
            \end{align}
            \vspace{-0.4cm} \normalsize
        \item Particle in-flight temperature \footnotesize
        \vspace{-0.2cm}
            \begin{align}
                \begin{split}
                f(\text{temperature} | \text{PFR}, \text{SOD}, \lambda, \text{CV}, \text{TGF}) &=     \exp\Big( \beta_0 
                    + \beta_1 \cdot \text{PFR} 
                    + \beta_2 \cdot \text{SOD} 
                    + \beta_3 \cdot \lambda 
                    + \beta_4 \cdot \text{TGF} \\
                    &+ \beta_5 \cdot \lambda^2 
                    + \beta_6 \cdot \text{TGF}^2 
                    + \beta_{7} \cdot \text{PFR} \cdot \text{TGF} \\
                    &+ \beta_{8} \cdot \text{SOD} \cdot \lambda
                    + \beta_{9} \cdot \text{SOD} \cdot \text{TGF} \Big)
            \end{split}
            \end{align}
            \vspace{-0.4cm} \normalsize
    \end{itemize}
    \item Process performance properties
    \begin{itemize}
        \item Deposition rate \footnotesize
        \vspace{-0.2cm}
            \begin{align}
            \begin{split}
                f(\text{rate} | \text{PFR}, \text{SOD}, \lambda, \text{CV}, \text{TGF}) &=     \exp\Big( \beta_0 
                    + \beta_1 \cdot \text{PFR} 
                    + \beta_2 \cdot \text{SOD} 
                    + \beta_3 \cdot \lambda 
                    + \beta_4 \cdot \text{TGF} \\
                    &+ \beta_5 \cdot \text{PFR}^2 
                    + \beta_6 \cdot \lambda^2 
                    + \beta_7 \cdot \text{CV}^2 
                    + \beta_{8} \cdot \text{TGF}^2 \\
                    &+ \beta_{9} \cdot \text{PFR} \cdot \text{TGF} \Big)
            \end{split}
            \end{align}
            \vspace{-0.4cm} \normalsize
        \item Deposition efficiency \footnotesize
        \vspace{-0.2cm} 
            \begin{align}
                \begin{split}
                f(\text{efficiency} | \text{PFR}, \text{SOD}, \lambda, \text{CV}, \text{TGF}) &=     \exp\Big( \beta_0 
                    + \beta_1 \cdot \text{PFR} 
                    + \beta_2 \cdot \text{SOD} 
                    + \beta_3 \cdot \lambda 
                    + \beta_4 \cdot \text{TGF} \\
                    &+ \beta_5 \cdot \lambda^2 
                    + \beta_6 \cdot \text{CV}^2 
                    + \beta_{7} \cdot \text{TGF}^2 
                    + \beta_{8} \cdot \text{PFR} \cdot \text{TGF} \Big)
            \end{split}
            \end{align}
            \vspace{-0.4cm} \normalsize
    \end{itemize}
    \item Coating quality properties
    \begin{itemize}
        \item Coating thickness \footnotesize
        \vspace{-0.2cm}
            \begin{align}
            \begin{split}
                f(\text{thickness} | \text{PFR}, \text{SOD}, \lambda, \text{CV}, \text{TGF}) &=     \exp\Big( \beta_0 
                    + \beta_1 \cdot \text{PFR} 
                    + \beta_2 \cdot \lambda 
                    + \beta_3 \cdot \text{CV} 
                    + \beta_4 \cdot \text{TGF} \\
                    &+ \beta_5 \cdot \lambda^2 
                    + \beta_6 \cdot \text{CV}^2 
                    + \beta_{7} \cdot \text{TGF}^2 
                    + \beta_{8} \cdot \text{PFR} \cdot \lambda \\
                    &+ \beta_{9} \cdot \text{CV} \cdot \text{TGF} \Big)
            \end{split}
            \end{align}
            \vspace{-0.4cm} \normalsize
        \item Coating roughness \footnotesize
        \vspace{-0.2cm}
            \begin{align}
                \begin{split}
                f(\text{roughness} | \text{PFR}, \text{SOD}, \lambda, \text{CV}, \text{TGF}) &=     \exp\Big( \beta_0 
                    + \beta_1 \cdot \text{PFR} 
                    + \beta_2 \cdot \text{SOD} 
                    + \beta_3 \cdot \lambda 
                    + \beta_4 \cdot \text{CV} \\
                    &+ \beta_5 \cdot \text{TGF} 
                    + \beta_6 \cdot \text{CV}^2
                    + \beta_7 \cdot \text{TGF}^2 \\
                    & + \beta_{8} \cdot \text{SOD} \cdot \text{TGF} 
                    + \beta_{10} \cdot \lambda \cdot \text{TGF} \Big)
            \end{split}
            \end{align}
            \vspace{-0.4cm} \normalsize
        \item Surface hardness \footnotesize
        \vspace{-0.2cm}
            \begin{align}
                \begin{split}
                f(\text{hardness} | \text{PFR}, \text{SOD}, \lambda, \text{CV}, \text{TGF}) &=     \exp\Big( \beta_0 
                    + \beta_1 \cdot \text{PFR} 
                    + \beta_2 \cdot \text{SOD} 
                    + \beta_3 \cdot \lambda 
                    + \beta_4 \cdot \text{CV} \\
                    &+ \beta_5 \cdot \text{TGF} 
                    + \beta_{6} \cdot \lambda \cdot \text{CV}
                    + \beta_{7} \cdot \lambda \cdot \text{TGF} \Big)
            \end{split}
            \end{align}
            \vspace{-0.4cm} \normalsize
        \item Coating porosity \footnotesize
        \vspace{-0.2cm}
            \begin{align}
                \begin{split}
                f(\text{porosity} | \text{PFR}, \text{SOD}, \lambda, \text{CV}, \text{TGF}) &=     \exp\Big( \beta_0 
                    + \beta_1 \cdot \text{PFR} 
                    + \beta_2 \cdot \text{SOD} 
                    + \beta_3 \cdot \lambda 
                    + \beta_4 \cdot \text{CV} \\
                    &+ \beta_5 \cdot \text{TGF} 
                    + \beta_6 \cdot \text{CV}^2 
                    + \beta_7 \cdot \text{TGF}^2
                    + \beta_{8} \cdot \text{PFR} \cdot \lambda \\
                    &+ \beta_{9} \cdot \text{PFR} \cdot \text{CV}
                    + \beta_{10} \cdot \text{PFR} \cdot \text{TGF}
                    + \beta_{11} \cdot \lambda \cdot \text{CV} \\
                    & + \beta_{12} \cdot \lambda \cdot \text{TGF} \Big)
            \end{split}
            \end{align}
    \end{itemize}
\end{enumerate} \normalsize

The coefficients for each gamma regression model, detailing the influence of each input parameter on the corresponding HVOF coating property, are presented in Table \ref{tab:coefficients}. These coefficients, found in \cite{rannetbauer2024predictive}, represent the contribution of each process parameter (PFR, SOD, $\lambda$, CV, TGF) to the respective HVOF coating properties. 

\begin{table}[h!]
\centering
\caption{Estimated coefficients for the gamma regression models of HVOF coating properties \cite{rannetbauer2024predictive}.}
\label{tab:coefficients}
\resizebox{\textwidth}{!}{
\begin{tabular}{l|rrrrrrrrrrrrr}
\hline\hline
\textbf{Property} & $\beta_0$ & $\beta_1$ & $\beta_2$ & $\beta_3$ & $\beta_4$ & $\beta_5$ & $\beta_6$ & $\beta_7$ & $\beta_8$ & $\beta_9$ & $\beta_{10}$ & $\beta_{11}$ & $\beta_{12}$ \\ \hline

Velocity             & \texttt{6.1297}      &\texttt{-0.0056}      & \texttt{0.0494}     & \texttt{-0.0182}      & \texttt{0.0483}      & \texttt{-0.0304}     & \texttt{-0.0219}      & \texttt{-0.0087}      & \texttt{0.0068}       & \texttt{0.0142}      &         &    &       \\ 
Temperature          & \texttt{7.4491}      &\texttt{-0.0072}      & \texttt{-0.0195}     & \texttt{0.0254}      & \texttt{0.0513}      & \texttt{-0.0056}     & \texttt{-0.0123}      & \texttt{-0.0042}      & \texttt{0.0040}       & \texttt{0.0034}      &         &    &       \\ 
 Rate      & \texttt{3.6337}      &\texttt{0.2668}      & \texttt{-0.0207}     & \texttt{0.0636}      & \texttt{0.1259}      & \texttt{-0.0303}     & \texttt{-0.0259}      & \texttt{-0.0540}      & \texttt{-0.0524}       & \texttt{0.0184}      &         &    &       \\ 
 Efficiency&  \texttt{-0.4546}      & \texttt{0.0051}      & \texttt{-0.0207}     & \texttt{0.0636}     &\texttt{0.1259}      & \texttt{-0.0273}     & \texttt{-0.0553}      & \texttt{-0.0538}     & \texttt{0.0184}       &      &         &    &       \\ 
Thickness            & \texttt{4.8928}      &\texttt{0.2275}      & \texttt{0.0664}     & \texttt{-0.2658}      & \texttt{0.0376}      & \texttt{-0.0338}     & \texttt{0.0428}      & \texttt{-0.0492}      & \texttt{-0.0254}       & \texttt{-0.0331}      &         &    &       \\ 
Roughness            & \texttt{3.5241}      &\texttt{0.0229}      & \texttt{-0.0065}     & \texttt{-0.0342}      & \texttt{-0.0419}      & \texttt{-0.0979}     & \texttt{0.0325}      & \texttt{0.0164}      & \texttt{-0.0219}       & \texttt{-0.0242}      &   \texttt{-0.0665}      &    &       \\ 
Hardness             & \texttt{6.3520}      &\texttt{-0.0372}      & \texttt{-0.0345}     & \texttt{0.0025}      & \texttt{-0.0192}      & \texttt{0.1189}     & \texttt{-0.0216}      & \texttt{0.0248}      &        &       &         &    &       \\ 
Porosity             & \texttt{2.7056}      &\texttt{0.0046}      & \texttt{0.0146}     & \texttt{-0.0293}      & \texttt{0.0074}      & \texttt{-0.0462}     & \texttt{0.0363}      & \texttt{0.0134}      & \texttt{0.0242}       & \texttt{0.0294}      &   \texttt{-0.0150}      &  \texttt{-0.0366}  &   \texttt{-0.0233}    \\         \hline \hline
\end{tabular}}
\end{table}

\end{document}

%% file: header.tex
\usepackage{booktabs}
\usepackage{amsmath}
\usepackage{amsthm}
\usepackage{amssymb}
\usepackage{cite}
\usepackage{url}
\usepackage[multiple]{footmisc}
\usepackage{graphicx}
\usepackage{epstopdf}
\usepackage{chngcntr}
\usepackage{enumitem}
\usepackage{hhline}
\usepackage{array}
\usepackage{hyperref}
\usepackage{color}
\usepackage{array}
\usepackage{multirow}
\usepackage{algorithm,algorithmic}
\usepackage[acronym]{glossaries}

\hypersetup{colorlinks=false}

\if@twoside
\else 
\fi

\numberwithin{equation}{section}
\counterwithin{figure}{section}
\counterwithin{table}{section}

\theoremstyle{plain}

\theoremstyle{definition}
\newtheorem{definition}{Definition}[section]

\newtheorem{example}{Example}[section]

\theoremstyle{remark}
\newtheorem{remark}{Remark}[section]

\definecolor{aog}{rgb}{0.0, 0.5, 0.0}

\makeatletter
\newcommand{\ALOOP}[1]{\ALC@it\algorithmicloop\ #1%
  \begin{ALC@loop}}
\newcommand{\ENDALOOP}{\end{ALC@loop}\ALC@it\algorithmicendloop}

\makeatother
